\documentclass[11pt]{article}
\usepackage{amssymb}
\usepackage{amsmath}
\usepackage{amsfonts}
\usepackage{amsthm}

\numberwithin{equation}{section}

\newtheorem{theorem}{Theorem}[section]
\newtheorem{lem}{Lemma}[section]
\newtheorem{prop}{Proposition}[section]
\newtheorem{corollary}{Corollary}[section]
\newtheorem{definition}{Definition}[section]
\newtheorem{remark}{Remark}[section]

\def\R{\mathbb{R}}
\def\Z{\mathbb{Z}}

\def\N{\mathbb{N}}

\def\T{\mathbb{T}}

\def\supp{\mathop{\rm supp}\nolimits}

\def\into{\int \hspace*{-4mm} - \,}

\begin{document}
\title{Sharp ill-posedness   results for the KdV and mKdV equations on the torus}
\author{\normalsize \bf  Luc Molinet  \\
{\footnotesize \it L.M.P.T., Universit\'e Fran\c cois Rabelais Tours, F\'ed\'eration Denis Poisson-CNRS,} \\
{\footnotesize \it Parc Grandmont, 37200 Tours,  France} }
\maketitle
{\bf AMS Subject Classification :} 35Q53,  35A01, 35D30. \vspace*{2mm}\\
\begin{abstract}
We establish a new a priori bound for $ L^2 $-bounded sequences of solutions to the mKdV equations on the torus.  
This first enables us to construct weak solutions in $ L^2$ for this equation and to check that the "solutions" constructed by Kappeler and Topalov in the defocusing case satisfy the equation in some weak sense. In a second time, we   prove that the solution-map associated with the mKdV and the KdV equation are discontinuous for the $ H^s(\T) $ topology for respectively $ s<0$ and $ s<-1$. These last results are sharp.
  \end{abstract}
\section{Introduction}
In this paper we study different properties of the Cauchy problems posed on the flat torus   $\T:=\R/2\pi \Z$ associated with the Korteweg-de Vries (KdV)  equation 
   \begin{equation}\label{KdV}
  w_t+w_{xxx} -6 w w_x =0 \; 
  \end{equation}
 and  the modified  Korteweg-de Vries  (mKdV) equation 
    \begin{equation}\label{mKdV}
    v_t+v_{xxx} \mp  6 v^2 v_x =0
  \end{equation}
   Here, $ w $ and $v $ are real-valued functions on $\T$. For some results we will have to distinguish between two mKdV equations depending on the sign  
    in front of the nonlinear term. \eqref{mKdV} is called  the defocussing mKdV equation  when there is a minus  sign in front of the nonlinear term and  the focussing mKdV equation when it is a plus sign.  The Cauchy problem associated with these equations in space of rough functions on the torus has been extensively studied these last two decades. In a seminal paper \cite{Bourgain1993}, Bourgain 
    proved that the Cauchy problem associated with the KdV equation is globally well-posed in $ H^s(\T) $, $s\ge 0 $, whereas the one associated with the mKdV equation is globally well-posed in $ H^s(\T) $, $ s\ge 1$,  and locally well-posed in $ H^s(\T) $ for $ s\ge 1/2$. The local well-posedness of the KdV equation was pushed down to $ H^s(\T) $, $s\ge -1/2 $ by Kenig, Ponce and Vega \cite{KPVJAMS96} (see \cite{CKSTT1} for the global-wellposedness 
     of the KdV and the mKdV equations in $ H^s(\T) $ for respectively $ s\ge -1/2 $ and $s\ge 1/2$.) The local well-posednesss results proved in these papers mean the following : for any initial data  $ u_0\in H^s(\T) $ there exists a time 
     $ T=T(\|u_0\|_{H^s})>0 $ only depending on $\|u_0\|_{H^s} $ and a solution  $ u$ that satisfies the equation  at least in  some weak sense and is unique in some function space $ X\hookrightarrow  C([0,T];H^s(\T))$. Moreover, for any $ R>0 $,  the flow-map $ u_0\mapsto u $ is continuous from the ball centered at the origin with radius $ R $ of $ H^s(\T) $ into $C([0,T(R)];H^s(\T)) $. 
     Note that in all these works, by a change of variables, the study of the KdV equation is actually restricted to initial data with mean value zero and  the mKdV equation is substituting by the following  "renormalized" equation :
     \begin{equation}\label{mKdV*}
       u_t+u_{xxx} \mp  6 (u^2-\into u^2) u_x =0
     \end{equation}
     where $ \into u^2 $ denotes the mean value of $ u^2 $. 
  The best  results quoted above are known to be sharp if one requires moreover the smoothness of the flow-map (cf.  \cite{Bo4})  or  the uniform continuity on bounded sets of the solution-map (cf. \cite{CCT2003}) Êassociated respectively with the KdV equation on space of functions with mean value zero and  with \eqref{mKdV*}. 
  On the other hand,  they  have been improved if one only requires the continuity of the flow-map. 
     In this direction, in \cite{KT1}-\cite{KT2}, Kappeler and Topalov introduced the following notion of solutions which a priori does not always corresponds to   the solution in the sense of distribution : 
     {\it ¬ A continuous curve $ \gamma \, :\, (a,b) \to H^\beta(\T) $ with $ 0\in (a,b)$ and $ \gamma(0)=u_0 $ is called a solution of KdV equation  (resp. mKdV equation) in $ H^\beta(\T) $ with initial data $ u_0 $ iff for any $ C^\infty $-sequence of initial data $ \{u_{0,n}\} $ converging to $ u_0 $ in $ H^\beta(\T) $ and for any $ t\in ]a,b[ $, the sequence of emanating solutions  $ \{ u_n\} $ of the KdV equation (resp. mKdV equation) satisfies : $ u_n(t) \to \gamma(t) $ in $H^\beta(\T) $¬}
   
     Note that a solution in the sense of this definition is necessarily unique. With this notion of solution they proved the global well-posedness of the KdV and the defocusing mKV equations in $ H^s(\T) $ for respectively $ s\ge -1 $ and 
     $ s\ge 0 $, with a solution-map which is continuous from $ H^{-1}(\T) $ (resp. $L^2(\T)) $ into $C(\R;H^{-1}(\T)) $ (resp. $C(\R;L^2(\T)) $).  Their proof is based on the inverse scattering method and thus depends in a crucial way of the complete integrability of these equations. It is worth noticing that, by Sobolev embedding theorem, their solutions of the defocussing mKdV equation satisfy the equation  in the  distributional sense as soon as $ s\ge 1/6$. 
     Independently, Takaoka and Tsutsumi (\cite{TT}) extended the local well-posedness of the mKdV equation (with the classical notion of solutions) 
      to $ H^s(\T) $ for $ s>3/8$ by modifying in a suitable way the Bourgain's space used as resolution space. This approach has been very recently improved by Nakanashi, Takaoka  and Tsutsumi \cite{NTT} and local well-posedness has been pushed to $ H^s(\T) $ for $ s>1/3 $ (local existence of solutions is shown in $ H^s(\T) $ for $ s>1/4$)
\vspace{2mm} 

In this paper we first establish an a priori estimate for $ L^2 $-bounded sequences of solutions to the mKdV equation. To this aim we slightly modify the spaces introduced by Ionescu-Kenig  and Tataru in \cite{IKT}. Recall that these spaces  are constructed by  localizing in time the Bourgain spaces with a localization in time that depends inversely on the space frequencies of the functions (see \cite{KT} and \cite{CCT} for previous works in this direction). Note that, to some extent, this approach is a  version  for the  Bourgain's spaces of the approach developped by Koch and Tzvetkov \cite{KTz} in Strichartz spaces. Once our a priori estimate  is established we translate it in the Bourgain's type spaces introduced in \cite{IKT} . This enables us to pass to the limit on the nonlinear term by separating resonant and non resonant parts. Following some ideas developped in \cite{M1}, we then derive a non continuity result for the mKdV equation in $ H^s(\T) $ for $ s<0 $. On the other hand, we obtain the existence of weak  $ L^2$-solutions of \eqref{mKdV*} and prove that the $ L^2$-solutions constructed in \cite{KT1} of the defocusing mKdV equation  satisfy the equation in some weak  sense. Finally, we follow some ideas of \cite{M2} and use properties of the Riccati map proved in \cite{KT0} to derive a non continuity result for the KdV equation in $ H^s(\T) $ for $ s<-1 $.  
 \subsection{Statement of the results}
      Our results can be summarized in the two following theorems. The first one deals with the discontinuity of the solution-map associated with the KdV and mKdV equations.
        \begin{theorem}\label{main}
        The Cauchy problems associated with the KdV equation and the mKdV equation are ill-posed in $ H^{s}(\T) $ for respectively $ s<-1 $ and $ s<0$. More precisely, 
        \begin{enumerate}
        \item[i)] for any $ T>0 $ and any $ s<-1 $, the solution-map $ u_0 \mapsto u $ associated with the KdV equation is  discontinous  at any $ u_0\in H^\infty_0(\T) $ 
        from $ H^\infty_0(\T) $, endowed with the topology inducted by $ H^s(\T) $, into  $ {\mathcal D}'(]0,T[\times \T) $. 
         \item[ii)]  for any $ T>0 $ and any $ s<0 $, the solution-map $ u_0 \mapsto u $ associated with the mKdV equation is  discontinous  at any non constant function  $ u_0\in H^\infty(\T) $    from $ H^\infty(\T) $, endowed with the topology inducted by 
         $ H^s(\T) $, into  $ {\mathcal D}'(]0,T[\times \T) $. 
\end{enumerate}
          \end{theorem}
\begin{remark} \label{rere}
Actually we prove the following assertions :
   \begin{enumerate}
        \item[i')] For any $ T>0 $,  the solution-map $ u_0 \mapsto u $ associated with the KdV equation is  discontinous, at any $ u_0\in H^\infty_0(\T) $, 
        from $ H^\infty_0(\T) $ endowed with the weak topology of  $ H^{-1}(\T) $ into  $ {\mathcal D}'(]0,T[\times \T) $. 
         \item[ii')]  For any $ T>0 $, the solution-map $ u_0 \mapsto u $ associated with the mKdV equation is  discontinous, at any non constant function  $ u_0\in H^\infty(\T) $,    from $ H^\infty(\T) $ endowed with the weak topology of 
         $ L^2(\T) $ into  $ {\mathcal D}'(]0,T[\times \T) $.
         \item[ii'')] Let $ u_0\in H^\infty(\T) $ be a  non constant function. There exists no $T>0 $ such that  for all $ t\in ]0,T[$ the flow-map $ u_0 \mapsto u(t) $ associated with the mKdV equation is  continuous, at   $ u_0\ $,    from $ H^\infty(\T) $ endowed with the weak topology of 
         $ L^2(\T) $ into  $ {\mathcal D}'(\T) $ .

\end{enumerate}
\end{remark} 
\begin{remark}
Our proof of the ill-posedness of the KdV equation below $ H^{-1}(\T) $ is heavily related to the algebraic structure of the equation via the Miura map. However, it is interesting to notice that, in \cite{MV}, a similar result is proved for the KdV-Burgers equation with a completely different method. 
\end{remark}
The second one deals with the existence of weak $L^2$-solutions to the mKdV equation.

           \begin{theorem}\label{Existence mKdV}\mbox{  } 
           \begin{enumerate}
 \item[i)] For any $ u_0\in L^2(\T) $ there exists a weak solution $ u\in C_{w}(\R;L^2(\T))\cap(\cup_{s<0}{\tilde F}^{s,1/2}) $ of  mKdV such that      $ u(0)=u_0 $. Moreover, $ u(t) \to u_0 $ in $ L^2(\T) $ as $ t\to 0 $.\\
  \item[ii)]  The $C(\R; L^2(\T))$-functions determined by the unique continuous extension to $ C(\R;L^2(\T)) $ of the $C(\R; H^\infty(\T))$ solution-map of the defocusing mKdV equation, constructed in \cite{KT2}, 	are weak  solutions of the mKdV equation and belong  to $ \cup_{s<0}{\tilde F}^{s,1/2}$. 
  \end{enumerate}
  \end{theorem}
  \begin{remark}
  The function spaces $ {\tilde F}^{s,1/2}  $ are defined by   \eqref{defF} when substituting the linear KdV group $ U(\cdot) $ by the linear group $ V(\cdot) $ defined in \eqref{taz}. Our notion of weak solution to mKdV is  
  described  in Definition \ref{def2}  (see also \eqref{mKdV3*} for assertion {\it ii))}.
  \end{remark}
 \begin{remark}\label{remark1}
  Once, the second assertion  of Theorem \ref{Existence mKdV}  is established, the first assertion seems to  have no more interest for the defocusing mKdV equation (note that  the first assertion  is, up to our knowledge, the only available existence result  of global  weak  $ L^2(\T) $-solutions  for the focusing mKdV equation). However, the proof of assertion 2 uses the complete integrability of the equation which is not a priori conserved by perturbations. On the other hand, the proof of assertion 1 seems to be widely more tractable and for instance certainly works for a wide class of perturbations of the defocusing mKdV equation. For instance,  
 $$
 \mbox{ Damped mKdV   :}\quad  u_t + u_{xxx}Ê+\nu u \mp  u^2 u_x =0 , \quad \nu>0  
 $$
 $$
\mbox{  KdV-mKdV : }\quad  u_t + u_{xxx}Ê\mp  u^2 u_x \mp  u u_x =0 ,
 $$
 \end{remark}
 \begin{remark} We  also construct in Proposition \ref{propo} global weak solutions for the "renormalized" mKdV equation \eqref{mKdV*}.
 \end{remark}
 {\it This paper is organized as follows:} In the next section we introduce the notations and the functions spaces  we will work with. We also give some useful estimates for time-localized functions. In Section 3 we recall general linear estimates in such functions spaces and some linear and bilinear estimates relating to the KdV group.  Section 4 is devoted to the proof of the uniform bound for $ L^2(\T) $-bounded sequence of solutions to mKdV. We prove Theorem 
  \ref{Existence mKdV} in Section 5 and Theorem \ref{main} in Section 6. Finally, in the appendix we first give a simplified proof of the continuous embedding  
  in $ L^4(\R\times \T) $ of some Bourgain's space related to the KdV group. Then, for sake of completeness, we prove some needed bilinear estimates and sketch the proof of some properties of the Riccati map  $ u\mapsto u'+u^2- \into u^2 $.
  \section{Notations and functional spaces} 
 \subsection{Notations}
 For $x,y\in \R_+^* $,  $x\sim y$ means that there exist $C_1$,
$C_2>0$ such that $C_1 x\leq y \leq C_2 x$.  $x\lesssim y$ and $ x\gtrsim y $ 
mean that there exists $C_2>0$ such that respectively  $ x \leq C_2 y $ and $x\ge C_2 y $. For
a Banach space $ X $, we denote by $ \| \cdot \|_X $ the norm in
$ X $.\\
We will use the same notations as in \cite{CKSTT1} and \cite{CKSTT2} to deal with Fourier transform of space periodic functions with a large period $\lambda $.
 $ (d\xi)_\lambda $ will be the renormalized counting measure on $ \lambda^{-1} \Z $ :
$$
\int a(\xi)\, (d\xi)_\lambda = \frac{1}{\lambda} \sum_{\xi\in
\lambda^{-1}\Z } a(\xi) \quad .
$$
As written in \cite{CKSTT2}, $ (d\xi)_\lambda $ is the counting measure on the integers when $ \lambda=1 $ and converges weakly
 to the Lebesgue measure when $ \lambda\to \infty $. In all the text, all the Lebesgue norms in $ \xi $ will be with respect to the
 measure $ (d\xi)_\lambda $. For  a $ (2\pi \lambda) $-periodic function $ \varphi$, we define its space Fourier transform on
  $ \lambda^{-1}\Z$ by
$$
\hat{\varphi}(\xi)=\int_{\lambda \T} e^{-i \xi x} \, f(x)
\, dx , \quad \forall \xi \in \lambda^{-1}\Z \quad .
$$
We denote by $ U(\cdot) $ the free group associated with the linearized Korteweg-de Vries equation,
$$
\widehat{U(t)\varphi}(\xi)=e^{i p(\xi) t} \,
\hat{\varphi}(\xi) , \quad \xi\in \lambda^{-1}\Z \quad , \; p(\xi)=\xi^3  \; .
$$
The Lebesgue spaces $ L^q_\lambda $, $ 1\le q \le \infty $,  for $ (2\pi\lambda)$-periodic functions, will be defined as usually by
$$
\|\varphi\|_{L^q_\lambda}=\Bigl( \int_{\lambda \T} |\varphi(x) |^q \, dx\Bigr)^{1/q}
$$
with the obvious modification for $ q=\infty $. \\
We define the Sobolev spaces $ H^s_{\lambda} $ for $ (2\pi\lambda)$-periodic functions by
$$
\|\varphi\|_{H^s_\lambda}=\|\langle \xi \rangle^{s} \varphi(\xi)
\|_{L^2_\xi} =\|J^s_x \varphi \|_{L^2_\lambda} \quad ,
$$
where $ \langle \cdot \rangle = (1+|\cdot|^2)^{1/2} $ and $
\widehat{J^s_x \varphi}(\xi)=\langle \xi \rangle^{s}
\widehat{\varphi}(\xi)
$. \\
Note that the closed subspace of zero mean value functions of $ H^s_\lambda $ will be denoted by $ {H}^s_{0,\lambda} $ (it is equipped with the $ H^s_\lambda $-norm).\\
In the same way, for a function $ u(t,x) $ on $\R\times \lambda\T $, we define its space-time Fourier transform by
$$
\hat{u}(\tau,\xi)={\cal F}_{t,x}(u)(\tau,\xi)=\int_{\R} \int_{\lambda \T} e^{-i (\tau
t+ \xi x)} \, u(t,x) \, dx dt , \quad \forall (\tau,\xi) \in
\R\times \lambda^{-1}\Z \quad .
$$
 $ L^p_t L^q_\lambda $ and $ L^p_T L^q_\lambda $ will denote respectively the Lebesgue spaces
$$
\|u\|_{L^p_t L^q_\lambda}=\Bigl( \int_{\R}
\|u(t,\cdot)\|_{L^q_\lambda}^p \, dt \Bigr)^{1/p}\quad \mbox{ and
}\quad \|u\|_{L^p_T L^q_\lambda}=\Bigr( \int_{0}^T
\|u(t,\cdot)\|_{L^q_\lambda}^p \, dt\Bigr)^{1/p}
$$
with the obvious modification for $ p=\infty $. \\
For any $ (s,b) \in \R^2 $, we define the  Bourgain space $ X^{s,b}_\lambda $,  of $ (2\pi\lambda) $-periodic (in $x$) functions 
  as the completion of $ {\mathcal S}(\lambda\T\times \R) $ for the norm 
\begin{eqnarray}
\| u \|_{X^{s,b}_{\lambda} }  =
 \| \langle \tau-p(\xi)\rangle^{b}  \langle \xi \rangle^s
  \hat{u}\|_{L^2_{\tau,\xi}} =
  \| \langle \tau\rangle^{b}  \langle \xi \rangle^s
  {\cal F}_{t,x}(U(-t) u ) \|_{L^2_{\tau,\xi}}\; ,
\end{eqnarray}
For $ T > 0$ and a function space $ B_\lambda$,
   we denote by  $ B_{T,\lambda} $ the corresponding restriction in time   space
 endowed with the norm
$$
\| u \|_{B_{T,\lambda}} =\inf_{w\in B_{\lambda}} \{ \|
w\|_{B_{\lambda}} , \, w(\cdot)\equiv u(\cdot) \hbox{ on } ]-T,T[ \, \}\;.
$$
Finally, for all function spaces of $ (2\pi\lambda) $-periodic functions, we will drop the index $ \lambda $ when $ \lambda=1 $. \\
 \subsection{Bourgain's spaces on frequency dependent time intervals}\label{section22}
 We will  need a Littlewood-Paley analysis. Let $\psi\in C^\infty_0(\R)$ be an even function  such that $\psi\geq 0$, $\supp \psi\subset [-3/2,3/2]$, $\psi\equiv 1$ on $[-5/4,5/4]$. We set $\eta_0:=\psi $ and  for all $k\in \N^* $, $\eta_k(\xi):=\psi(2^{-k} \xi)-\psi(2^{-k+1} \xi)$ and $ \eta_{\le k}:=\psi(2^{-k} \cdot)=\sum_{j=0}^k \eta_k $. We also set  $\tilde{\eta}_0:=\psi(2^{-1}\cdot)  $ and  for all $k\in \N^* $, $\tilde{\eta}_k(\xi):=\psi(2^{-k-1} \xi)-\psi(2^{-k+2} \xi)$
 \\
 The Fourier multiplicator operators by $ \eta_j $, $\tilde{\eta}_j $ and $ \eta_{\le j} $  will be denoted respectively by $ \Delta_j $, $\tilde{\Delta}_j$ and $ S_j$,  and the  projection on the constant Fourier mode will be denoted by $ P_0 $ , i.e. for any $ u\in L^2(\lambda\T)$
  $$
 \widehat{\Delta_j u}:= \eta_j \widehat{u} , \quad  \widehat{\tilde{\Delta}_j u}:= \tilde{\eta}_j \widehat{u},
 \quad \widehat{S_j u}:=\eta_{\le j} \hat{u} 
 \quad \mbox{Êand }Ê\quad P_0(u) =\frac{1}{2\pi \lambda} \int_{\lambda \T} u(x)\, dx  .
$$
By a slight abuse of notations, we will also define the operator $ \Delta_j  $, $\tilde{\Delta}_j$  and $ P_0 $ on $ L^2(\lambda\T \times \R)$-functions by 
 the same formula.  Finally, for any $ l\in \N $ we define the functions $ \nu_l $ on $ \lambda^{-1} \Z\times \R $ by 
 \begin{equation}\label{defnunu}
  \nu_l(\xi,\tau):=\eta(\tau-p(\xi)) \; . 
  \end{equation}
   Let $ 0\le b<1 $, $ k\in\N$, $ t_0\in \R $ and $ f\in C^\infty(\lambda \T\times ]t_0-2^{-k},t_0+2^{-k}[Ê) $, we define
 $$
 \| f\|_{ F_{\lambda,k, t_0}^b} := \inf_{{\tilde f}\in X^{0,b,1}_\lambda} \Big\{ \| {\tilde f} \|_{ X^{0,b,1}_\lambda} , \; {\tilde f}=f \mbox{ on }  ]t_0-2^{-k},t_0+2^{-k}[Ê\Bigr\} 
 $$
 where $ X^{0,b,1}_\lambda $ is the Bourgain's type space defined by 
 $$
  X^{0,b,1}_\lambda=\left\{ \begin{array}{l} 
 f\in {\mathcal S}'(\lambda\T\times \R)   ,\\
\displaystyle  \|Êf\|_{X^{0,b,1}_\lambda} := \sum_{j=0}^\infty 2^{j b } \| \nu_j(\xi,\tau) \widehat{f} \|_{L^2_{\xi,\tau}} < \infty 
 \end{array} \right\}\; .
 $$
  Our a priori estimate will take place in the normed space $ G_\lambda $ defined as the completion of  $  C^\infty(\lambda \T\times \R) $ for the norm
  \begin{equation}
 \|Êf\|_{G_\lambda}^2 :=\sup_{t\in \R} \sum_{k\ge 0}  \| \Delta_k f \|_{F_{\lambda,k,t}^{1/2}}^2 \; .  \label{defG}
 \end{equation}
Once our a priori estimate will be established we will make use of the spaces $ F_\lambda^{s,b} $, introduced in \cite{IKT}, that are endowed with the norm 
\begin{equation}\label{defF}
 \|Êf\|_{F_\lambda^{s,b}}^2 := \sum_{k\ge 0} \Bigl( \sup_{t\in \R} 2^{ks}  \| \Delta_k f \|_{F_{\lambda,k,t}^b}\Bigr)^2 \; .
\end{equation}
To handle the nonlinear term, for $ k\in \N $, we will also need to introduce  the function space $ Z^b_{\lambda,k} $
 defined as the completion of $ L^2(\lambda\T\times \R) $ for the following norm :
 \begin{equation}
 \| f\|_{Z^b_{\lambda,k}} =\sum_{j\ge 0} 2^{bj}Ê\| \nu_j(\xi,\tau) \langle \tau-p(\xi)+i 2^k \rangle^{-1} \widehat{f} \|_{L^2} \; .
 \end{equation}
 Finally, for  $ t_0\in \R $ and $ f\in L^2(\T\times ]t_0-2^{-k},t_0+2^{-k}[Ê) $, we define
 $$
 \| f\|_{ Z^b_{\lambda,k,t_0}} := \inf_{{\tilde f}\in Z^b_{\lambda,k}} \Big\{ \| {\tilde f} \|_{Z^b_{\lambda,k}} , \; {\tilde f}=f \mbox{ on }  ]t_0-2^{-k},t_0+2^{-k}[Ê\Bigr\} \; .
 $$
 \subsection{Some useful estimates for localized in time functions }
 Following  \cite{IKT} , we state the lemma below that  will be useful in the  linear estimates and also in the nonlinear estimates when we will localize the functions on time interval of length of order $ 2^{-l} $. 
 \begin{lem}
 Let be given $ b\in [0,1[ $, $\lambda\ge 1$ and $ f\in X^{0,b,1}_\lambda $. Then for all $ l\in \N $ it holds 
 \begin{equation}\label{ds1} 
 2^{bl} \Bigl\| \eta_{\le l}(\tau-p(\xi))  \int_{\R} |\widehat{f}(\xi,\tau')| \, (2^{-l} (1+2^{-l}|\tau-\tau'|)^{-4} \, d\tau'    \Bigr\|_{L^2_{\xi,\tau}} 
 \lesssim   (1\vee 2^{ (b-1/2)l} )\|f\|_{X^{0,b,1}_\lambda}
 \end{equation}
 and
 \begin{align}
 \sum_{j\ge l+1} 2^{bj} \Bigl\| \eta_{j }(\tau-p(\xi))   \int_{\R} |\widehat{f}(\xi,\tau')| &  \, 2^{-l} (1+2^{-l}|\tau-\tau'|)^{-4} \, d\tau'    \Bigr\|_{L^2_{\xi,\tau}} \nonumber \\
  & \lesssim   (1\vee 2^{ (b-1/2)l} ) \|f\|_{X^{0,b,1}_\lambda} \; .\label{ds2} 
 \end{align}
 \end{lem}
 \proof  For $ b\in [0,1/2] $  and $ l\in \N $, we get by Cauchy-Schwarz in $ \tau' $,
 \arraycolsep1pt
\begin{eqnarray*}
I^{b,l} & := & 
 2^{bl} \Bigl\| \eta_{\le l}(\tau-p(\xi))  \int_{\R} |\widehat{f}(\xi,\tau')| \, 2^{-l} (1+2^{-l}|\tau-\tau'|)^{-4} \, d\tau'    \Bigr\|_{L^2_{\xi,\tau}} \nonumber\\
 & \lesssim &  2^{bl} \Bigl\| \eta_{\le l}(\tau-p(\xi)) \sum_{q=0}^\infty \int_{\R} 
  \nu_q(\xi,\tau')|\widehat{f}(\xi,\tau')| \, 2^{-l} (1+2^{-l}|\tau-\tau'|)^{-4} \, d\tau'    \Bigr\|_{L^2_{\xi,\tau}}  \nonumber\\
  & \lesssim & 2^{(b-1)l} \Bigl\| \eta_{\le l}(\tau-p(\xi)) \sum_{q=0}^\infty  I^{b,l}_q
\Bigl\| \nu_q(\xi,\tau') |\widehat{f}(\xi,\tau')| \langle \tau'-p(\xi)\rangle^{b}  \Bigr\|_{L^2_{\tau'}}    \Bigr\|_{L^2_{\xi,\tau}} 
\end{eqnarray*}
\arraycolsep4pt
where
\begin{eqnarray*}
I_q^{b,l} & := & \sup_{|\tau-p(\xi)| \le 2^{l+1}} \Bigl( 
\int_{\R} {\tilde{\eta}_q(\tau'-p(\xi))} (1+2^{-l}|\tau-\tau'|)^{-8} \langle \tau'-p(\xi)\rangle^{-2b} \, d\tau'  \Bigr)^{1/2}
\end{eqnarray*}
But, for $ q\le l+2 $, 
$$
I_q^{b,l} \lesssim   \Bigl( \int_{\R} {\tilde{\eta}_q(\tau'-p(\xi))} \langle \tau'-p(\xi)\rangle^{-2b} \, d\tau'  \Bigr)^{1/2}
\lesssim (2^q 2^{-2bq})^{1/2} \lesssim 2^{(\frac{1}{2}-b)l}  
$$
and for $ q\ge l+3 $, noticing that $|\tau-\tau'|Ê\sim 2^q >\! >2^l $ in the region where the integrand is not vanishing, we obtain
$$
I_q^{b,l} \lesssim (2^q 2^{-2bq} 2^{8(l-q)})^{1/2} \lesssim 2^{(\frac{1}{2}-b)l}\; .
$$
Gathering the above estimates, we eventually get
$$
I^{b,l}\lesssim 2^{(b-1)l} 2^{l/2}Ê 2^{(\frac{1}{2}-b)l} \| f\|_{X^{0,b,1}_\lambda} =\| f\|_{X^{0,b,1}_\lambda}\; ,
$$
which proves \eqref{ds1}Ê for $ b\in [0,1/2]Ê$. The case $ b>1/2 $ follows immediately by using the result for $ b=1/2 $ and controlling the 
 $ X^{0,1/2,1}_\lambda$-norm by the  $ X^{0,b,1}_\lambda$-norm.

Now to prove \eqref{ds2}, we first notice that by the mean-value theorem, 
$$
|\eta_j(\tau-p(\xi))-\eta_j (\tau'-p(\xi))|\lesssim  2^{-j} |\tau-\tau'| \
$$
and it thus  would be sufficient to  estimate
$$
\sum_{j\ge l+1} 2^{bj} \Bigl\|  \int_{\R}   \eta_{j }(\tau'-p(\xi)) |\widehat{f}(\xi,\tau')| \, 2^{-l} (1+2^{-l}|\tau-\tau'|)^{-4} \, d\tau'    \Bigr\|_{L^2_{\xi,\tau}} 
$$
$$
+\sum_{j\ge l+1} 2^{(b-1)j} \Bigl\|  \int_{\R} {\tilde \eta}_{j }(\tau'-p(\xi))  |\widehat{f}(\xi,\tau')| \, 2^{-l} \frac{|\tau-\tau'|}{(1+2^{-l}|\tau-\tau'|)^{4}} \, d\tau'    \Bigr\|_{L^2_{\xi,\tau}} \; .
$$
By identifying a convolution term and applying generalized Young's inequality, we can bound the first term for any $ b\ge 0  $ by 
$$
 \sum_{j\ge l+1} 2^{bj} \|Ê\nu_{j}  \widehat{f} \|_{L^2} \int_{\R} 2^{-l} (1+2^{-l}|\tau'|)^{-4} \, d\tau' \lesssim \sum_{j\ge l+1} 2^{bj} \|Ê\nu_{j}  \widehat{f} \|_{L^2}
$$
We call by $ J^{b,l} $ the second term. For $ b\in [0,1/2] $ we proceed as for \eqref{ds1} to get 
$$
J^{b,l}\lesssim 2^{-l} \sum_{j\ge l+1} 2^{(b-1)j} \Bigl\| {\tilde \eta}_{j}(\tau-p(\xi)) \sum_{q=0}^\infty  J_{q,j}^{b,l}
\Bigl\| \nu_q(\xi,\tau') |\widehat{f}(\xi,\tau')| \langle \tau'-p(\xi)\rangle^{b}  \Bigr\|_{L^2_{\tau'}} \Bigr\|_{L^2_{\tau,\xi}}
$$
with
 \begin{eqnarray*}
J_{q,j}^{b,l} & := & \sup_{\langle \tau-p(\xi)\rangle \in [ 2^{j-2}, 2^{j+2}]} \Bigl( 
\int_{\R} {\tilde{\eta}_q(\tau'-p(\xi))} |\tau-\tau'|^2 (1+2^{-l}|\tau-\tau'|)^{-8} \langle \tau'-p(\xi)\rangle^{-2b} \, d\tau'  \Bigr)^{1/2}
\end{eqnarray*}
For $ q\le j-5$, we use that $ |\tau-\tau'|\sim 2^j \ge 2^l $ in the region where the integrand is not vanishing  to  get 
\begin{eqnarray*}
J_{q,j}^{b,l} & \lesssim & 2^{j} 2^{4(l-j)}  \Bigl( \int_{\R} {\tilde{\eta}_q(\tau'-p(\xi))} \langle \tau'-p(\xi)\rangle^{-2b} \, d\tau'  \Bigr)^{1/2} \\
& \lesssim &  2^{j} 2^{4(l-j)} 2^{(\frac{1}{2}-b)q} \lesssim 2^{4l} 2^{-(\frac{5}{2}+b)j}
\end{eqnarray*}
whereas for $ q\ge j-4 $, we easily get  
$$
J_{q,j}^{b,l} \lesssim  2^{-bj}\Bigr\|\frac{|\tau|}{(1+2^{-l}|\tau|)^{4}}\Bigl\|_{L^2_\tau} Ê\lesssim 2^{-bj} 2^{\frac{3l}{2}} .
$$
Gathering these estimates we conclude that 
$$
J^{b,l}\lesssim \|f\|_{X^{0,b,1}_\lambda}  2^{-l} 2^{\frac{3l}{2}}  \sum_{j=l+1}^\infty2^{(b-1)j} 2^{j/2}  2^{-bj} 
\lesssim  \|f\|_{X^{0,b,1}_\lambda}  2^{l/2}   \sum_{j=l+1}^\infty 2^{-j/2} \lesssim \|f\|_{X^{0,b,1}_\lambda}  \; .
$$
 This proves \eqref{ds2} for $ b\in [0,1/2]$. Finally, in the case  $ b\in ]1/2,1]$  we observe that 
 $$
J^{b,l}\lesssim 2^{-l} \sum_{j\ge l+1} 2^{(b-\frac{1}{2})j} 2^{-j/2} \Bigl\| {\tilde \eta}_{j}(\tau-p(\xi)) \sum_{q=0}^\infty  
 J_{q,j}^{1/2,l}
\Bigl\| \nu_q(\xi,\tau') |\widehat{f}(\xi,\tau')| \langle \tau'-p(\xi)\rangle^{1/2}  \Bigr\|_{L^2_{\tau'}} \Bigr\|_{L^2_{\tau,\xi}}
$$
  and the above bounds on $ J_{q,j}^{1/2,l} $ lead to 
 $$
J^{b,l}
\lesssim  \|f\|_{X^{0,1/2,1}_\lambda}  2^{l/2}   \sum_{j=l+1}^\infty 2^{(b-\frac{1}{2})j } 2^{-j/2} \lesssim 2^{(b-\frac{1}{2})l}  \|f\|_{X^{0,b,1}_\lambda}  
  .
  $$ \qed
\begin{remark} A very useful corollary of the preceding lemma is the following:  Let  $ b\in [0,1/2] $, $ f\in X^{0,b,1}_\lambda $ and $ \gamma\in C^\infty_c(\R) $ with support in $ ]-2,2[ $. Then  for all $ k\in \N $, it holds 
\begin{equation}\label{gh1}
 \|\gamma(2^k t) f\|_{X^{0,b,1}_\lambda} \lesssim  \|Êf\|_{X^{0,b,1}_\lambda} 
\end{equation}
and 
\begin{equation}\label{gh2}
\| \eta_{\le k} (\tau-p(\xi)) {\mathcal F}_{xt}(\gamma(2^k t) f)\|_{L^2_{\xi,\tau}} \lesssim 2^{-bk } \|Êf\|_{X^{0,b,1}_\lambda} \; .
\end{equation}
\end{remark} 
\begin{remark} \label{remark2}
It is easy to check that $ G_\lambda $ is continuously embedded in $ L^\infty(\R; L^2_\lambda(\T)) $. Indeed, for any $ u\in G_\lambda$, $ t_0\in\R 
$ and $ k\in \N $, taking a function  $ {\tilde u}\in X^{0,1/2,1}_\lambda $ such that $ {\tilde u}\equiv u $ on $ ]t_0-2^{-k},t_0+2^{-k}[ $ and 
$$
\|{\tilde u} \|_{X^{0,1/2,1}_\lambda} \le  2 \| u\|_{F_{\lambda,k,t_0}^{1/2}} \; , 
$$
it holds
$$
{\mathcal F}_x (\Delta_k u(t_0))(\xi)=\int_{\R}  {\mathcal F}( {\Delta_k \tilde u}(\xi,\tau) ) e^{it_0 \tau}Ê\, d\tau 
$$
According to   the obvious estimate 
\begin{equation}\label{ds3}
\Bigl\| \int_{\R} |\widehat{f}(\xi,\tau') |\, d\tau' \Bigr\|_{L^2_\xi} \lesssim \| f\|_{X^{0,1/2,1}_\lambda} , \quad \forall f\in X^{0,1/2,1}_\lambda, 
\end{equation}
 this leads to 
$$
\|\Delta_k  u (t_0)\|_{L^2_\lambda} \lesssim \| \Delta_k {\tilde u}   \|_{X^{0,1/2,1}_\lambda}\lesssim \|\Delta_k u\|_{F^{1/2}_{\lambda,k,t_0}}\; .
$$
Squaring and summing in $ k $ one obtains that
\begin{equation}\label{lll}
 \|u\|_{L^\infty_t L^2(\lambda \T))} \lesssim \|u\|_{G_\lambda} \; .
 \end{equation}
On the other hand, it seems pretty clear that $ G_\lambda $ is not included in $C(\R; L^2_\lambda(\T)) $.
\end{remark} 
 \section{Some  linear and bilinear estimates} 
 \subsection{General linear estimates}
  We first derive linear estimates that do not depend on the 
  dispersive linear group associated with  our functional space. We mainly follow \cite{KT}-\cite{IKT}. 
 \begin{lem}\label{linear1}
 Let be given $ b\in [1/2,1[$. Then $\forall \varphi\in L^2_\lambda $ and all $ k\in \N $, it holds 
 \begin{equation}
 \| U(t) \varphi \|_{F^{b}_{\lambda,k,0}} \lesssim 2^{(b-\frac{1}{2})k} \| \varphi\|_{L^2_\lambda} \; .
 \end{equation}
 \end{lem}
 \proof Clearly, it suffices to prove that for any $ k\in \N $,
  $$
  \| \eta_0(2^k  t) U(\cdot)  \varphi \|_{X^{0,b,1}_\lambda}\lesssim 2^{(b-\frac{1}{2})k}  \| \varphi\|_{L^2_\lambda} \; .
 $$
 Notice that the left-hand side member of the above inequality is bounded  by 
 $$
\Bigl( \sum_{j=0}^\infty 2^{bj} \| \eta_j(\tau) 2^{-k}Ê\hat{\eta}_0(2^{-k}\tau) \|_{L^2_{\tau}} \Bigr)\|  \varphi\|_{L^2_\lambda} \; .
 $$
 Since $ \eta_0\in C^\infty_c(R) $, $ \widehat{\eta_0} $  decays at least as $ (1+|y|)^{-4} $ and thus 
 \begin{equation}\label{gsgs}
  \|  \eta_j(\tau) 2^{-k}Ê\hat{\eta}_0(2^{-k}\tau)\|_{L^2_{\tau}}\lesssim  
  \|  \eta_j(\tau) 2^{-k}Ê(1+2^{-k} |\tau|)^{-4})\|_{L^2_{\tau}}
   \lesssim 2^{-k} 2^{j/2} \min (1, 2^{4(k-j)}) \; .
 \end{equation}
 Hence, 
 $$
 \sum_{j=0}^{k+2} 2^{bj} \|Ê\eta_j(\tau) 2^{-k}Ê\hat{\eta}_0(2^{-k}\tau)  \|_{L^2_\tau}ÊÊ\lesssim 2^{(b-\frac{1}{2})k} 
 $$
 and 
 $$
 \sum_{j=k+2}^{\infty} 2^{bj} \|Ê\eta_j(\tau) 2^{-k}Ê\hat{\eta}_0(2^{-k}\tau) \|_{L^2}Ê\lesssim  
 \sum_{j=k+2}^{\infty} 2^{j(b-\frac{7}{2})} 2^{3k}Ê\lesssim 2^{(b-\frac{1}{2})k}  \; . \qed
 $$
 \begin{lem}\label{linear2}
 Let be given $ b\in [1/2,1[$. Then  for any $ k\in \N $ and any $ f\in Z^b_{\lambda,k,0} $ it holds 
 \begin{equation}
 \Bigl\|\int_0^t U(t-t') f (t') \, dt' \Bigr\|_{F^{b}_{\lambda,k,0}} \lesssim 2^{(b-\frac{1}{2})k}\|  f \|_{Z^b_{\lambda,k,0}}\; .
 \end{equation}
 \end{lem}
 \proof  Let $ {\tilde f }\in Z^b_{\lambda,k}  $ be and extension of $ f$ such that 
 $ \| {\tilde f}\|_{Z^b_{\lambda,k}} \le 2 \| f\|_{Z^b_{\lambda,k,0}} $. 
 We set 
 $$
 v:={\eta_0}(2^k t) \int_0^t U(t-t')  {\tilde f}  (t')\, dt' \; .
 $$ 
 Then 
 \begin{eqnarray*}
 {\cal F}_{tx} (v) (\xi,\tau) & = & {\cal F}_{t}\Bigl[ \eta_0(2^k t) \int_{\R} \frac{e^{it\tau}-e^{itp(\xi)}}{i(\tau-p(\xi))} {\cal F}_{t,x}({\tilde f}) \, d\tau \Bigr] (\xi,\tau)\\
 &= &    {\cal F}_{t}(\eta_0(2^k t) )\ast \Bigl[  \frac{{\cal F}_{t,x}({\tilde f})}{i(\tau-p(\xi))} \Bigr]-
 {\cal F}_{t}(e^{itp(\xi)}\eta_0(2^k t) )\int_{\R} \frac{{\cal F}_{t,x}({\tilde f})}{i(\tau'-p(\xi))}\, d\tau' \\
 & =  &\int_{\R} {\cal F}_{t,x}({\tilde f})(\tau',\xi) \Bigl[ \frac{2^{-k} \hat{\eta}_0 (2^{-k}(\tau-\tau'))-2^{-k} \hat{\eta}_0 (2^{-k} (\tau-p(\xi)))} {i(\tau'-p(\xi))}\Bigr] \, d\tau' 
 \end{eqnarray*}
 Now we claim that $$
 I=\frac{\Bigl|2^{-k}\hat{\eta}_0(2^{-k}(\tau-\tau'))-2^{-k}\hat{\eta}_0(2^{-k}(\tau-p(\xi)))\Bigr|}{|\tau'-p(\xi)|}|\tau'-p(\xi)+i2^k|
 $$
 \begin{equation}\label{tg}
 \lesssim 2^{-k} (1+2^{-k}|\tau-\tau'|)^{-4}+ 2^{-k} (1+2^{-k}|\tau-p(\xi)|)^{-4}
 \end{equation}
 Assuming \eqref{tg} for a while, it follows that  
 \begin{eqnarray*}
 \|v\|_{X^{0,b,1}} & \lesssim  & \sum_{j=0}^\infty 2^{bj}Ê\Bigl\| \eta_j(\tau-p(\xi)) 
 \int_{\R} \frac{\widehat{{\tilde f}}(\xi,\tau') }{\langle \tau'-p(\xi) +i 2^k \rangle}
 2^{-k}  (1+2^{-k}|\tau-\tau'|)^{-4} \, d\tau' \Bigr\|_{L^2_{\xi,\tau} }\\
 & & + \sum_{j=0}^\infty 2^{bj}Ê
 Ê\Bigl\| \eta_j(\tau-p(\xi))  2^{-k}  (1+2^{-k}|\tau-p(\xi)|)^{-4}
 \int_{\R} \frac{\widehat{{\tilde f}}(\xi,\tau') }{\langle \tau'-p(\xi) +i 2^k \rangle}\, d\tau'
  \Bigr\|_{L^2_{\xi,\tau}}
 \end{eqnarray*} 
 The desired bound on the first term of the above right-hand side member follows directly from  \eqref{ds1}-\eqref{ds2}. To obtain the desired bound on the second term, we combine \eqref{ds3} and \eqref{gsgs}.\\
  It thus remains to prove \eqref{tg}. For this we first notice that,  since  $ \eta_0(y)=\eta_0(|y|) $, by the mean-value theorem there exists $ \theta\in ]|\tau-\tau'|,|\tau-p(\xi)|[$ such that 
 \begin{equation}\label{tg2}
 \Bigl|2^{-k}\hat{\eta}_0(2^{-k}(\tau-\tau'))-2^{-k}\hat{\eta}_0(2^{-k}(\tau-p(\xi)))\Bigr|\lesssim 
 2^{-2k} \hat{\eta}_0'(2^{-k}\theta) |\tau'-p(\xi)| \; .
 \end{equation}
Furthermore, since $ \eta_0 \in {\mathcal  S}(\R) $, 
$$
 |  \hat{\eta}_0(y)|+ |  \hat{\eta}_0'(y)|\lesssim (1+|y|)^{-10} \; .
 $$
  Let us now  separate three cases : 
  \begin{enumerate}
\item[$\bullet$] $ |\tau-p(\xi)|\le  2^k $. Then \eqref{tg}  is obvious whenever $ |\tau'-p(\xi)|\ge 2^k $ and follows directly from \eqref{tg2} whenever  $ |\tau'-p(\xi)|\le 2^k $ . 
\item[$\bullet$]  $ |\tau-p(\xi)|\ge  2^k $ and $|\tau-\tau'|\sim |\tau-p(\xi)|   $. Then  we must  have $ |\theta|\sim |\tau-p(\xi)| $ and 
 $ |\tau'-p(\xi)+i 2^k|\lesssim |\tau-p(\xi)| $.Therefore   \eqref{tg2} leads to  
$$
I\lesssim 2^{-2k} (1+2^{-k} |\tau-p(\xi)|)^{-5} |\tau-p(\xi)|Ê\lesssim 2^{-k} (1+2^{-k} |\tau-p(\xi)|)^{-4}Ê
$$
\item[$\bullet$] $ |\tau-p(\xi)|\ge  2^k $ and  $|\tau-\tau'|\not \sim |\tau-p(\xi)|   $ .   Then $ |\tau'-p(\xi)|\sim (|\tau-p(\xi)|\vee |\tau'-\tau|) \gtrsim 2^k $ and \eqref{tg} Êfollows directly from the decay of $\hat{\eta}_0 $.\qed
\end{enumerate}

\subsection{Specific linear and bilinear estimates}
We will also need estimates that are specific for Bourgain's spaces associated with the KdV linear group. We first recall the following Strichartz's type estimate proved in \cite{Bourgain1993} (we give a simplified proof of this estimate in the appendix) : 
\begin{lem}\label{lem1}
For any  $ \lambda\ge 1 $ and  any $ u\in X^{0,1/3}_\lambda $, it holds 
\begin{equation}\label{strichartz}
\|u\|_{L^4_{t,\lambda}} \lesssim \|u\|_{X^{0,1/3}_\lambda}
\end{equation}
\end{lem}
Finally we will  make a frequent use of  the following bilinear estimates that can be deduced for instance from \cite{Tao}  (we give a  proof of these estimates in the appendix since we need to quantify the dependence of these estimates with respect to the period $ \lambda $): 
\begin{lem} \label{lem2}
Let $ \lambda\ge 1 $ and let $ u_1 $ and $ u_2 $ be two real valued $ L^2 $ functions defined on $ \R\times ( \lambda^{-1}\Z) $ with the following support properties
\begin{equation*}
(\tau,\xi)\in \supp u_i \Rightarrow 
\langle \tau-\xi^3 \rangle \lesssim L_i , \, i=1,2.
\end{equation*}
Then for any $ N>0  $  the following estimates holds:
\begin{equation} \label{propro1}
\| u_1\star u_2\|_{L^2_\tau L^2(|\xi|\ge N)} \lesssim(L_1\wedge L_2)^{1/2}  \Bigl(  \frac{(L_1\vee  L_2)^{1/4}}{N^{1/4}} +\lambda^{-1/2}\Bigr)  \| u_1\|_{L^2}  \|  u_2\|_{L^2} \; ,
\end{equation}
and
\begin{equation} 
 \| \Lambda[N](u_1,u_2) \|_{L^2_{\tau,\xi}}
 \lesssim (L_1\wedge  L_2)^{1/2} \Bigl(  \frac{(L_1\vee L_2)^{1/2}}{N} +\lambda^{-1/2}\Bigr)  \| u_1\|_{L^2}  \|  u_2\|_{L^2}. \label{propro2}
\end{equation}
where $ \Lambda[N] : (L^2(\R\times \lambda^{-1}\Z))^2 \to L^\infty(\R\times \lambda^{-1}\Z) $ is defined by 
$$
\Lambda[N](u_1,u_2)(\tau,\xi)=\int_{\R} \int_{\Bigl||\xi_1|-|\xi-\xi_1|\Bigr|\ge \frac{N}{100}}u_1(\tau_1,\xi_1) u_2(\tau-\tau_1,\xi-\xi_1)\, (d\xi_1)_\lambda\, d\tau_1  
$$
\end{lem}
\section{ A priori estimate for     smooth solutions to \eqref{mKdV*}}
\label{4.1}
As in    previous works on mKdV on the torus (cf. \cite{Bourgain1993},\cite{CKSTT1}), we actually  work with the "renormalized " mKdV equation \eqref{mKdV*} instead of the mKdV equation itself. This permits to cancel some  resonant part in  the nonlinear term.  Recall that for $ v \in C(\R; H^\infty_\lambda)$,  a smooth solution to mKdV with initial data $ v_0$,  the $ L^2 $-norm of $ v$ is a constant of the motion  and thus   $ u(t,x):=v(t,x\mp \frac{6t}{2\pi \lambda} \|v_0\|_{L^2_\lambda}^2 ) $ satisfies \eqref{mKdV*}.

Denoting by $N(u) $ the  nonlinear term of \eqref{mKdV*}, it holds for any $ \xi \in \lambda^{-1}\Z $, 
\begin{eqnarray*}
{\mathcal F}_{x} [N(u)](\xi)& = & -\frac{6i}{\lambda^2} \sum_{\xi_1+\xi_2+\xi_3=\xi \atop \xi_1+\xi_2\neq 0 } \hat{u}(\xi_1)  \hat{u}(\xi_2) 
\xi_3  \hat{u}(\xi_3) \nonumber \\
& = & - \frac{2 i \xi}{\lambda^2}  \Bigl[ \sum_{\xi_1+\xi_2+\xi_3=\xi \atop (\xi_1+\xi_2)(\xi_1+\xi_3)(\xi_2+\xi_3) \neq 0}   \hat{u}(\xi_1)  \hat{u}(\xi_2)  \hat{u}(\xi_3) \Bigr]
\nonumber \\
& & + \frac{6  i \xi}{\lambda^2}  \hat{u}(\xi)  \hat{u}(\xi)  \hat{u}(-\xi) \nonumber \\
& := &  -i\xi \Bigl[  \Bigl( {\mathcal F}_{x}\Bigl[ A(u,u,u)\Bigr](\xi) + {\mathcal F}_{x} \Bigr[B(u,u,u)\Bigr] (\xi) \Bigr)\Bigr]  \, 
\end{eqnarray*}
i.e. 
\begin{equation}\label{defAB}
6(u^2-P_0(u^2))u_x=\partial_x \Bigl( A(u,u,u)+B(u,u,u) \Bigr) \; .
\end{equation}
According to the resonance relation \eqref{resonance}, $ A$ is non resonant whereas $ B $ is a resonant term. As pointing out  in \cite{TT}, $ A$ is a "good term" as far as one wants  to solve the equation in $ H^s(\T) $ for $ s\ge 1/4 $. On the other hand,  $ B $ is  a bad term as soon as one wants to solve the equation below $ H^{1/2} (\T)$, giving rise to rapid oscillations that breaks the uniform continuity on bounded set of the flow-map. 
\begin{prop}\label{prolo}
Let $ \lambda\ge 1 $ and  $ u\in C(\R; H^\infty(\lambda \T)) $ be a solution to \eqref{mKdV*}. Then, 
\begin{equation}\label{estF}
\|u\|_{G_\lambda}\lesssim \| u\|_{L^\infty(\R; L^2_\lambda)} + \|Êu\|_{G_\lambda}^3 
\end{equation}
\end{prop}Ê
 \proof 
 From the definition of the norm $ G_\lambda $, we have to bound $ \displaystyle \sup_{t_0\in \R} \sum_{k\ge 0}\|\Delta_k u \|_{F^{1/2}_{\lambda,k,t_0}}^2 $.
   We use that for any $( t_0,t)  \in \R^2 $, it holds 
 $$
  u(t)= U(t-t_0)  u(t_0)+\int_{t_0}^t U(t-t') \partial_x \Bigl(A({u}(t'))+B({u}(t')\Bigr)dt' \; .
 $$
 By translation in time we can always assume that $ t_0=0 $ and 
 according to Lemmas \ref{linear1}-\ref{linear2},  
 $$
 \|U(t)\Delta_k u(0)\|Ê_{F^{1/2}_{\lambda,k,0}}Ê\lesssim \| \Delta_k u (0)\|_{L^2_\lambda}
 $$
 and 
 $$
 \Bigl\|\int_{0}^t U(t-t') \partial_x \Delta_k (A(u(t'))+B(u(t'))dt' \Bigr\|_{F^{1/2}_{\lambda,k,0}} \lesssim
  \| \Delta_k \partial_x (A(u))+B(u)) \|_{Z^{1/2}_{\lambda,k,0}}\; .
 $$
 Since $ \sum_{k\ge 0 } \| \Delta_k u (0)\|_{L^2_\lambda}^2\sim \| u\|_{L^\infty(\R; L^2_\lambda)}^2 $, it remains to prove that 
 $$
 \sum_{k=0}^\infty  2^{2k}\Bigl(  \| \Delta_k  A({u})\|_{Z^{1/2}_{\lambda, k,0}}^2+
 \|  \Delta_k B(u) \|_{Z^{1/2}_{\lambda, k,0}}^2\Bigr)  \lesssim  \| u\|_{G_\lambda}^6\; .
 $$
 This is the aim of the two following lemmas. 
 \begin{lem} \label{bil2}
 Let  be given $t_0\in \R$, $\lambda\ge 1 $ and $ u_i\in G_\lambda  $ for $ i=1,2,3$. Then it holds 
 $$
 \sum_{k=0}^\infty 2^{2k}  \| \Delta_k \Bigl(B({u}_1,{u}_2,{u}_3)\Bigr) \|_{Z^{1/2}_{\lambda,k,t_0}}^2 \lesssim  \prod_{i=1}^3 \| u_i\|_{G_\lambda}^2
 \; .
 $$
 \end{lem}
 \proof By translation in time, we can take $ t_0=0 $. 
 For any  fixed $ k\in \N $,  we take a time  extension $ {\tilde u_1} $ of $ u_1 $ such that $ \|\Delta_k {\tilde u_1}\|_{X^{0,1/2,1}_\lambda} \le 2 \|Ê\Delta_k u_1 \|_{F^{1/2}_{\lambda,k,0}}$. Then, in view of the structure of $ B $,  it holds 
 $$
 2^k \| \Delta_k \Bigl(B({u}_1,{u}_2,{u}_3)\Bigr) \|_{Z^{1/2}_{\lambda,k,0}}
 \lesssim \sum_{l\ge 0} 2^{l/2} 2^k \Bigl\| \eta_l (\tau-\xi^3)\langle \tau-\xi^3+i 2^k \rangle^{-1}      
 {\mathcal F}_{xt} \Bigl[\Delta_k B(v_1,v_2,v_3)\Bigr]\Bigr\|_{L^2} 
 $$
 $$
 \lesssim  \sum_{l\ge 0} 2^{l/2} 2^{k} \Bigl\| \eta_l (\tau-\xi^3)\langle \tau-\xi^3+i 2^k \rangle^{-1}  
  {\mathcal F}_{xt} \Bigl[  B(\Delta_k v_1,{\tilde \Delta}_k v_2,{\tilde \Delta}_k v_3)]\Bigr\|_{L^2} 
$$
where $ v_1 =\eta_0(2^k t ) {\tilde u_1} $ and $ v_i = \eta_0(2^k t) u_i  $, $ i=2,3 $. By duality it suffices to prove that  
\begin{eqnarray}\label{estB}
I_k& := & 2^k \Bigl| \Bigl(  {\mathcal F}_{xt} \Bigl[ B(\Delta_k v_1,{\tilde \Delta}_k v_2,{\tilde \Delta}_k v_3) \Bigr],
 \langle \tau-\xi^3+i 2^k \rangle^{-1}     \widehat{w} \Bigr)_{L^2}
   \Bigr| \nonumber \\
   & \lesssim&   \sup_{l} (2^{-l/2} \| \nu_l \widehat{w}\|_{L^2})  \|\Delta_k v_1 \|_{X^{0,1/2,1}_\lambda} 
     \prod_{i=2}^3 \|S_{k+1} v_i\|_{X^{0,1/2,1}_\lambda}  \; .\label{zb}
  \end{eqnarray} 
     Indeed, first, according to \eqref{gh1} for any $ k\in \N^* $ and any $u\in C(\R\times \lambda\T) $, 
     \begin{equation}\label{rg1}
     \| \eta_0(2^kt) S_{k-1} u\|_{X_\lambda^{0,1/2,1}}^2 \sim \sum_{j=0}^{k-1} \| \eta_0(2^kt)   \Delta_j u \|^2_{X^{0,1/2,1}_\lambda}  
     \le \sum_{j=0}^{k-1}Ê\|Ê\Delta_j u \|_{F^{1/2}_{\lambda,j,0}}^2 \le \|u\|_{G_\lambda}^2 \; .
     \end{equation}
    Second,  taking a $ C^\infty $-function $ \gamma \, :\, \R  \mapsto [0,1] $ with compact support in $ [-1,1] $ satisfying 
    $ \gamma\equiv 1$  on $ [-1/2,1/2] $ and $ \sum_{m\in \Z} \gamma(t-m)\equiv 1 $ on $\R$,  we get  
  for any $ j\ge k  \in \N  $ ,  \begin{eqnarray}
  \| \eta_0(2^k t)  \Delta_j u\|_{X^{0,1/2,1}_\lambda} &\le & \sum_{|m|\le 2^{j-k}} \|\gamma(2^{j}t -m) \eta_0(2^k t) \Delta_j u\|_{X^{0,1/2,1}_\lambda} \nonumber\\
& \le & \sum_{|m|\le  2^{j-k}}  \| \Delta_j u\|_{F^{1/2}_{\lambda,j,2^{-j} m }}\lesssim  2^{j-k}  \|u\|_{G_\lambda}
\; . \label{rg2}
\end{eqnarray}   
   Therefore,  \eqref{zb} will lead to 
    $$
    2^k \| \Delta_k \Bigl(B({u}_1,{u}_2,{u}_3)\Bigr) \|_{Z^{1/2}_{\lambda,k,0}}\lesssim \|\Delta_k u_1 \|_{X_{\lambda,k,0}} \prod_{i=2}^3  \| u_i \|_{G_\lambda} 
    $$
    which will gives the result by squaring and summing in $ k$.

          Since the norms in the right-hand side of \eqref{zb}Ê only see the size of the modulus of the 
          Fourier transform of the functions we can assume that all the functions have non negative Fourier transforms. 
          In view of  the structure of $ B $, using  Cauchy-Schwarz, we get 
              $$
     I_k   \lesssim  2^k  \| \eta_k  {\widehat v}_1 \ast {\tilde \eta}_k  {\widehat v}_2 \|_{L^2(\langle\xi\rangle \sim 2^k)} \| \Bigl(\langle \tau-\xi^3+i 2^k \rangle^{-1}  {\tilde \eta}_k \widehat{w} \Bigr) \ast {\tilde \eta}_k  \check{\widehat v}_3 \|_{L^2(\langle\xi\rangle\sim 2^k)} 
     $$
     where $ \check{\widehat{v}}(\tau,\xi)=\widehat{v}(-\tau,-\xi)$  for all $ (\tau,\xi)\in \R\times \lambda^{-1}\Z$.
     For $ k=0,1,2 $ this yields directly the result by using the Strichartz inequality \eqref{strichartz}. For $ k\ge 3 $ we introduce the following  notation :  we set 
     \begin{equation}\label{notations}
      {\tilde \nu}_l(\tau,\xi) :=\eta_l(\tau-\xi^3) \mbox{ for } l> k \; \mbox{ and } \quad  {\tilde \nu}_l(\tau,\xi) :=\sum_{j=0}^k \eta_j(\tau-\xi^3) \mbox{  for } l=k\; . 
      \end{equation}
      According to Lemma \ref{lem2}Ê we obtain 
     \begin{eqnarray}
 I_k    & \lesssim   & 2^k\sum_{\min(l_1,l_2,l_3) \ge k} \sum_{l\ge 0}Ê  2^{(l_1\wedge l_2)/2}   \Bigl(2^{(l_1\vee l_2)/4 }2^{-k/4} +1\Bigr) \|  {\tilde \nu}_{l_1} \eta_k  \widehat{v}_1\|_{L^2}  \|  {\tilde \nu}_{l_2}  {\tilde \eta}_k  \widehat{v}_2\|_{L^2}
  \nonumber \\
  & &  2^{(l\wedge l_3)/2}   \Bigl(2^{(l\vee l_3)/4 }2^{-k/4} +1\Bigr)2^{-(l\vee k)}Ê \|   {\tilde \nu}_{l_3}  {\tilde \eta}_k   \widehat{v}_3\|_{L^2}  \| \nu_{l}  {\tilde \eta}_k  \widehat{w}\|_{L^2} \nonumber \\
   & \lesssim   &  2^k \sum_{l\ge 0}  2^{\frac{3l}{8}} 2^{-(l\vee k)}Ê  \| \nu_{l}  {\tilde \eta}_k  \widehat{w}\|_{L^2} 
   \sum_{\min(l_1,l_2,l_3) \ge k}  \|  {\tilde \nu}_{l_1} \eta_k  \widehat{v}_1\|_{L^2}  \|  {\tilde \nu}_{l_2}  {\tilde \eta}_k  \widehat{v}_2\|_{L^2}
Ê \|   {\tilde \nu}_{l_3}  {\tilde \eta}_k   \widehat{v}_3\|_{L^2} \nonumber \\
& & \max\Bigl( 2^{\frac{3}{8}(l_1+l_2+l_3)} 2^{-k/2}, 2^{\frac{(l_1+l_2)}{4} } 2^{\frac{3l_3}{8}} 2^{-k/4},
    2^{\frac{3}{8}(l_1+l_2)} 2^{-k/4} 2^{\frac{l_3}{8}}, 2^{\frac{(l_1+l_2)}{4}}  2^{\frac{l_3}{8}} \Bigr)    \nonumber \\
  & \lesssim  &  \sup_{l}(2^{-l/2}  \| \nu_{l}  {\tilde \eta}_k   \widehat{w}\|_{L^2} ) \sup_{l\ge k }(2^{l/2}  \|  {\tilde \nu}_{l }  \eta_k  
   \widehat{v}_1\|_{L^2} )  \prod_{i=2}^3 \sup_{l\ge k }Ê\Bigl(  2^{l/2} \|  {\tilde \nu}_{l }  {\tilde \eta}_k   \widehat{v}_i\|_{L^2}\Bigr) \nonumber \\
 & \lesssim &  \sup_{l}(2^{-l/2}  \| \nu_{l }  {\tilde \eta}_k   \widehat{w}\|_{L^2} )
   \|\Delta_k  v_1 \|_{X^{0,1/2,1}_\lambda}  \prod_{i=2}^3 \|{S}_{k+1} v_i\|_{X^{0,1/2,1}_\lambda} \; \label{kkk}
     \end{eqnarray}
     where we used \eqref{gh2} in the last step.   \qed
      \begin{lem} \label{estimA}
      \label{bil1}Let be given $ t_0\in \R$, $ \lambda\ge 1 $ and   $ u_i \in  G_\lambda$ for $ i=1,2,3$. Then it holds 
 \begin{equation}
J:= \sum_{k=0}^\infty 2^{2k}  \| \Delta_k  \Bigl(A({u}_1,{u}_2,{u}_3)\Bigr) \|_{Z^{1/2}_{\lambda,k,t_0}}^2 \lesssim  \prod_{i=1}^3 \| u_i\|_{G_\lambda}^2\; .
 \label{estA}
 \end{equation}
 \end{lem}
 \proof Again by translation in time, we can take $ t_0=0$. 
 Denoting  by $ \xi_i  $ the Fourier modes of $ u_i $, we can always assume  by symmetry that $ |\xi_1|\ge |\xi_2|\ge |\xi_3 |$. 
 We divide $ A $ in different terms corresponding to regions of  $ (\lambda^{-1} \Z)^3 $. \\
{\bf 1.}Ê{$ |\xi_1|\le 2^4$}.\\
 Then it holds $ |\xi_1+\xi_2+\xi_3 |\le   2^6$.    By Sobolev inequalities  and \eqref{strichartz}, 
 $$
 J \lesssim \sum_{k=0}^7 \Bigl[ \prod_{i=1}^3 \| \eta_0(2^k t) S_6 u_i\|_{X^{0,1/2,1}_\lambda} \Bigr]^2
 \; , 
 $$
 which  is acceptable thanks to \eqref{rg1}-\eqref{rg2}.\vspace*{2mm} \\
 {\bf 2.}Ê$|\xi_1|\ge 2^4$ and   $ |\xi_1|\le 4  |\xi|$.\\
 In this region, it holds $ |\xi|\sim |\xi_1| $.  
  We rewrite $ \eta_k(\xi) $ as $ \eta_k(\xi_1) + \eta_k(\xi) -\eta_k(\xi_1) $ and notice that by the mean-value theorem, 
  \begin{equation}
 | \eta_k(\xi) -\eta_k(\xi_1) | \lesssim 
 \min\Bigl(1, 2^{-k}  \Bigl| |\xi|-|\xi_1|\Bigr| \Bigr)  \; . \label{tt}
 \end{equation}
  Therefore, $ J \lesssim\displaystyle  \sum_{k=0}^\infty  (J_{1,k}^2 +  J_{2,k}^2 )$ with 
  $$
  J_{1,k}:= 2^k \Bigl\|{\tilde \eta}_k(\xi)   \langle \tau-\xi^3+i2^{k}\rangle^{-1}  
 {\mathcal F}_{tx} \Bigl( A_{1} (\Delta_k u_1,u_2,u_3)\Bigr)  \Bigr\|_{F^{1/2}_{\lambda,k,0}}
  $$
  and 
  $$
  J_{2,k} := 2^k   \Bigl\|{\tilde \eta}_k(\xi) \langle \tau-\xi^3+i2^{k}\rangle^{-1}  
 {\mathcal F}_{tx} \Bigl( A_{2} (u_1,u_2,u_3)\Bigr)  \Bigr\|_{F^{1/2}_{\lambda,k,0}}\; ,
  $$
  where 
  \begin{equation}\label{defA1}
 {\mathcal F}_x\Bigl(  A_1(v_1,v_2,v_3)\Bigr) (\xi):=\frac{1}{\lambda^2} \sum_{(\xi_1,\xi_2,\xi_3)\in \Theta(\xi) \atop  2^4\le |\xi_1|\le 4 |\xi|} \hat{v}_1(\xi_1)  \hat{v}_2(\xi_2)  \hat{v}_3(\xi_3)\; ,
  \end{equation}
  $$
  {\mathcal F}_x\Bigl(  A_2(v_1,v_2,v_3)\Bigr) (\xi):=\frac{1}{\lambda^2} \sum_{(\xi_1,\xi_2,\xi_3)\in \Theta(\xi) \atop  2^4\le |\xi_1|\le 4 |\xi|}    [\eta_k(\xi)-\eta_k(\xi_1)]\ \hat{v}_1(\xi_1)  \hat{v}_2(\xi_2)  \hat{v}_3(\xi_3)\; ,
  $$
  and 
  $$
  \Theta(\xi):=\Bigl\{(\xi_1,\xi_2,\xi_3)\in \lambda^{-1}\Z ,\; \sum_{i=1}^3 \xi_i=\xi, \, \prod_{i,j=1\atop i\neq j}^3 (\xi_i+\xi_j)\neq 0  
    \mbox{ and  }  |\xi_3|\le |\xi_2|\le|\xi_1|\Bigr\} \; .
  $$
{\bf $ \bullet $ Estimate on $J_{1,k}$} \\
 For any  fixed $ k\in \N $,  we take a time  extension $ {\tilde u_1} $ of $ u_1 $ such that 
 $ \|\Delta_k {\tilde u_1}\|_{X^{0,1/2,1}_\lambda} \le 2 \|Ê\Delta_k u_1 \|_{F^{1/2}_{\lambda,k,0}} $ .
  We set $ v_1=\eta_0(2^k t) \tilde{u}_1 $ and  $ v_i=\eta_0(2^k t) u_i$, $ i=2,3 $.  By duality it suffices to prove that 
  \arraycolsep3pt
 \begin{eqnarray}
 H_{1,k}& := & 2^k \Bigl| \Bigl( {\mathcal F}_{xt} \Bigl[ A_1(\Delta_k v_1, v_2, v_3) \Bigr],
 \langle \tau-\xi^3+i 2^k \rangle^{-1}     \widehat{\tilde{\Delta} w} \Bigr)_{L^2}
   \Bigr| \nonumber \\
   & \lesssim&   \sup_{l} (2^{-l/2} \| \nu_l w\|_{L^2})  \|\Delta_k v_1 \|_{X^{0,1/2,1}_\lambda} 
     \prod_{i=2}^3   \|S_{k+3} v_i\|_{X^{0,1/2,1}_\lambda} \label{zzzb}
  \end{eqnarray} 
    which is acceptable thanks to  \eqref{rg1}-\eqref{rg2}.  Again, since the norms in the right-hand side of \eqref{zzzb}Ê only see the size of the modulus of the 
          Fourier transform of the functions we can assume that all the functions have non negative Fourier transforms. 
We separate four cases :\\
 {\bf A.} $\xi \xi_1 \le 0 $. Then $ |\xi_2+\xi_3|=|\xi-\xi_1|\sim |\xi| $. 
 Proceeding as in \eqref{kkk}, we get 
 \begin{eqnarray}
 H_{1,k} & \lesssim & 2^k \| \Bigl(\langle \tau-\xi^3+i 2^k \rangle^{-1} \tilde{\eta}_k \widehat{w}\Bigr) \ast {\eta}_{k}  \check{\widehat{v_1}} \|_{L^2(|\xi|\sim 2^k)} 
 \|Ê {\eta}_{\le k+3} \widehat{v_2} \ast  \eta_{\le k+3} \widehat{v_3}\|_{L^2(|\xi|\sim 2^k)}Ê
 \nonumber \\
& \lesssim &  \sup_{l}(2^{-l/2}  \| \nu_{l }  \widehat{{\tilde \Delta}_k  w}\|_{L^2} ) \|\Delta_k v_1 \|_{X^{0,1/2,1}_\lambda} 
 \prod_{i=2}^3     \|S_{k+3} v_i\|_{X^{0,1/2,1}_\lambda}\; .\end{eqnarray} 
  {\bf B.} $\xi \xi_2\le 0  $.  Then $ |\xi_1+\xi_3|=|\xi-\xi_2|\sim |\xi|$.  Therefore, exchanging the role of $ v_1 $ and $ v_2 $,  we can proceed exactly as in the previous case. \\
   {\bf C.} $\xi \xi_3\le 0  $.  Then $ |\xi_1+\xi_2|=|\xi-\xi_3|\sim |\xi|$.  Therefore, exchanging the role of $ v_1 $ and $ v_3 $,  we can proceed exactly as in the case {\bf A}. \\
 {\bf D.}  $ \xi_1,\, \xi_2 $ and $ \xi_3 $ are of the same sign. Then 
 $ |\xi-\xi_3|=|\xi_1+\xi_2|\sim |\xi| $. Therefore we can proceed exactly as in the previous case.\vspace{2mm} \\
   {\bf $ \bullet $ Estimate on $J_{2,k}$} \\
  For any fixed $ k\in \N $,  we set $ v_i:=\eta_0(2^k t) u_i $, $ i=1,2,3 $.
By duality and  \eqref{rg1}-\eqref{rg2}, it suffices to prove that 
 \begin{eqnarray}\label{estBB}
 H_{2,k}& := & 2^k \Bigl| \Bigl( {\mathcal F}_{xt} \Bigl[ A_2( v_1, v_2, v_3) \Bigr], 
 \langle \tau-\xi^3+i 2^k \rangle^{-1}    \widehat{ \tilde{\Delta}_k w} \Bigr)_{L^2}
   \Bigr| \nonumber \\
   & \lesssim&   2^{-k/8} \sup_{l} (2^{-l/2} \| \nu_l \widehat{w}\|_{L^2})  
      \prod_{i=1}^3     \|S_{k+3} v_i\|_{X^{0,1/2,1}_\lambda}  \label{zzb} \end{eqnarray} 
     Since the norms in the right-hand side of \eqref{zzb}Ê only see the size of the modulus of the 
          Fourier transform of the functions we can assume that all the functions have non negative Fourier transforms.  
    In view of \eqref{tt}, we thus infer that 
    $$  
       H_{2,k} \lesssim  2^k \Bigl| \Bigl(  {\mathcal F}_{xt} \Bigl[ \tilde{A_2}( v_1, v_2, v_3) \Bigr], 
 \langle \tau-\xi^3+i 2^k \rangle^{-1}        \widehat{\tilde{\Delta}_k w} \Bigr)_{L^2} \Bigr| \;  
   $$
   where 
  $$
  {\mathcal F}_x\Bigl(  \tilde{A_2}(v_1,v_2,v_3)\Bigr) (\xi):=\frac{1}{\lambda^2} \sum_{(\xi_1,\xi_2,\xi_3)\in \Theta(\xi) \atop  2^4\le |\xi_1|\le 4 |\xi|}   \Bigl[1\wedge (2^{-k}  \Big||\xi|-|\xi_1|\Bigr|\Bigr]   \hat{v}_1(\xi_1)  \hat{v}_2(\xi_2)  \hat{v}_3(\xi_3)\; .
  $$
First to estimate the contribution of the region  $  \Bigl| |\xi|-|\xi_1|\Bigr| \le |\xi|^{3/4}  $,
  we proceed exactly as for $ J_{1,k} $ by separating the four cases  {\bf A}, {\bf B}, {\bf C} and {\bf D} to  obtain 
 $$
 H_{2,k} \lesssim 2^{-k/4}  \sup_{l}(2^{-l/2}  \| \nu_{l }  \widehat{\tilde{ \Delta}_k w}\|_{L^2} )  \prod_{i=1}^3 
    \|S_{k+3} v_i\|_{X^{0,1/2,1}_\lambda}
$$
which is acceptable. \\
Now in the region $  \Bigl| |\xi|-|\xi_1|\Bigr| \ge |\xi|^{3/4}  $, we notice that  $ |\xi_2+\xi_3|=|\xi-\xi_1|\ge   |\xi|^{3/4}  $. Therefore, setting 
 \begin{equation}\label{defsigma}
 \sigma:=\sigma(\tau,\xi)=\tau-\xi^3 \mbox{ and  }\sigma_i:=\sigma(\tau_i,\xi_i), \; i=1,2,3,
 \end{equation}
 we claim that  the well-known resonance relation 
  \begin{equation}\label{resonance}
  \sigma-\sigma_1-\sigma_2-\sigma_3=3(\xi_1+\xi_2)(\xi_1+\xi_3)(\xi_2+\xi_3)
  \end{equation}
   leads to  (recall that $ |\xi|\sim |\xi_1|$)
\begin{equation}\label{ps}
\max(|\sigma|,|\sigma_i|) \gtrsim \lambda^{-1} |\xi|^{7/4}Ê \; 
\end{equation}
Indeed, either we are in the cases {\bf B} ,{\bf C}  or {\bf D} described above and then there exists $ i\in {2,3} $ such that 
 $ |\xi_1+\xi_i|\sim |\xi| $ so that we are done   or we are in the case {\bf A}. In this last case, we first notice   that $ |\xi_2+\xi_3|\sim |\xi| $. Second, since $ \xi  \xi_1 \le 0 $ and  $  \Bigl| |\xi|-|\xi_1|\Bigr| \ge |\xi|^{3/4}  $ we must have $ |\xi|\le |\xi_1| -|\xi|^{
3/4} $ and $ |\xi|=|\xi_2|+|\xi_3|-|\xi_1| $. It follows that  
  $ |\xi_3|\le |\xi_1|-\frac{1}{2} |\xi|^{3/4} $ and ensures that  \eqref{ps} holds also in this region. \\
 Using \eqref{propro1}-\eqref{propro2}
 and  the notations \eqref{notations}, we get 
 \begin{eqnarray*}
 H_{2,k} & \lesssim & 2^k \|Ê\eta_{\le k+2}   \widehat{v_2} \ast  \eta_{\le k+3}   \widehat{v_3}
 \|_{L^2(|\xi| \gtrsim 2^{3k/4} )}Ê \\
& & \Bigl\|\Lambda[2^{3k/4}] \Bigl(\langle \tau-\xi^3+i 2^k \rangle^{-1} {\tilde \eta}_k {\widehat{w}}\, , \,   \eta_{\le k+3} 
  \check{\widehat{v_1}}\Bigr) \Bigr\|_{L^2} \nonumber \\
   & \lesssim   & \sum_{\min(l_1,l_2,l_3) \ge k, l\ge 0 \atop \max(2^l,2^{l_i})\ge \lambda^{-1} 2^{7k/4}} Ê 2^k 2^{(l_2\wedge l_3)/2}   \Bigl(2^{(l_2\vee l_3)/4}2^{-\frac{3}{16}k} +\lambda^{-1/2}\Bigr)  \|  \tilde{\nu}_{l_2} \eta_{\le k+3}  \widehat{v_2}\|_{L^2} 
    \|  \tilde{\nu}_{l_3}  \eta_{\le k+3}  \widehat{v_3}\|_{L^2} \\
  & &  2^{(l_1\wedge l)/2}   
   \Bigl(2^{(l_1\vee l)/4 }2^{-\frac{3}{16}k} +\lambda^{-1/2}\Bigr)^{1/4} \Bigl(2^{(l_1\vee l)/2 }2^{-3k/4} +\lambda^{-1/2}\Bigr)^{3/4} 2^{-(l\vee k)}Ê\\
    & & \hspace*{20mm}  \|  \tilde{\nu}_{l_1}  \eta_{\le k+3}   \widehat{v_1}\|_{L^2}  \| \nu_{l}  {\tilde \eta}_k  \widehat{w}\|_{L^2}\\
  & \lesssim  & 2^{-\frac{k}{32}} \sup_{l}(2^{-l/2}  \| \nu_{l }  {\tilde \eta}_k  \widehat{w}\|_{L^2} ) \prod_{i=1}^3
   \sup_{l\ge k }Ê\Bigl(  2^{l/2} \|  \tilde{\nu}_{l }   S_{k+3}  v_i\|_{L^2_\lambda}\Bigr) \\
  & \lesssim &    2^{-\frac{k}{32}} \sup_{l}(2^{-l/2}  \| \nu_{l }  \widehat{{\tilde \Delta}_k  w}\|_{L^2} )  \prod_{i=1}^3 
   \|S_{k+3} v_i\|_{X^{0,1/2,1}_\lambda}
     \end{eqnarray*} 
     where in the last step we used \eqref{gh1}-\eqref{gh2}. \vspace{2mm} \\
{\bf 3.}Ê$|\xi_1|\ge 2^4 $ and $ |\xi_1|\ge 4 |\xi|$.\\
In this region it holds 
\begin{equation}\label{31}
|\xi_2+\xi_3|\ge \frac{3}{4} |\xi_1| , \; |\xi_2|\ge \frac{3}{8} |\xi_1| \mbox{ and  }  |\xi_3|\le \frac{5}{8}|\xi_1 |
\end{equation}
The two  first above inequalities are clear.   To prove the third one, we  notice that in this region $ \xi_1 \xi_2\le 0  $ and we proceed by contradiction by assuming that 
$ |\xi_3|> \frac{5}{8}|\xi_1 |$. Then we  first notice that  if $ \xi_1\xi_3 \ge 0 $  then we must have 
 $ |\xi|\ge |\xi_3| >\frac{5}{8}|\xi_1 |$ which contradicts $ |\xi_1|\ge 4 |\xi|$. Second, if $ \xi_1 \xi_3\le 0  $ then we have $ |\xi|=\Bigl||\xi_2|+|\xi_3|-|\xi_1|\Bigr|>\frac{5}{4}|\xi_1|-|\xi_1| =\frac{1}{4} |\xi_1| $ which again contradicts
$ |\xi_1|\ge 4  |\xi|$  .\\ 
Therefore  the resonance relation yields
\begin{equation}\label{32}
\max(|\sigma|,|\sigma_i|) \gtrsimÊ\lambda^{-1} |\xi_1|^2  \; .
\end{equation}
In this region $ \displaystyle J \lesssim \sum_{k\in \N} J_{3,k}^2 $ where, for any fixed $ k \in\N$, $ J_{3,k} $ is defined by 
\begin{equation}\label{defJ3}
J_{3,k}  :=Ê \sum_{k_1\ge k} 
\sum_{l\ge 0} 2^{l/2}  2^k \Bigl\| \eta_l (\tau-\xi^3)\langle \tau-\xi^3+i 2^k \rangle^{-1}  
 \tilde{\eta}_k(\xi){\mathcal F}_{tx} \Bigl(A_{3}(\Delta_{k_1} v_1,v_2,v_3)\Bigr) \Bigr\|_{L^2}
\end{equation}
with
\begin{equation}\label{defA3}
{\mathcal F}_x\Bigl(A_3(v_1,v_2,v_3)\Bigr)(\xi):= \frac{1}{\lambda^2}\sum_{(\xi_1,\xi_2,\xi_3)\in \Lambda(\xi) \atop   |\xi_1|\ge 2^4, \, 
 |\xi_1|\ge 4 |\xi|}   \hat{v}_1(\xi_1)  \hat{v}_2(\xi_2)  \hat{v}_3(\xi_3)
\end{equation}
and  $ v_i:=\eta_0(2^k t) u_i $, $ i=1,2,3 $. \\
Let $ \gamma $  be a $ C^\infty $-function $ \gamma \, :\, \R  \mapsto [0,1] $ with compact support in $ [-1,1] $ satisfying 
    \begin{equation}\label{defgamma}
     \gamma\equiv 1 \mbox{ on } [-1/2,1/2] \mbox{ and } \sum_{m\in \Z} \gamma(t-m)^3\equiv 1 \mbox{ on } \R.
     \end{equation}
     We set 
   $$ v_i^{k_1,m}: = \eta_0(2^k t)\gamma(2^{k_1}t -m) u_i =\gamma(2^{k_1}t -m) v_i, \;  i=1,2,3 , \, m\in \Z\; .
   $$
   Clearly, $ A_3(v_1,v_2,v_3)\equiv\displaystyle \sum_{|m|\lesssim 2^{k_1-k}} A_3( v_1^{k_1,m},v_2^{k_1,m},v_3^{k_1,m}) $. Therefore, by duality it suffices to prove that 
\begin{eqnarray}\label{estAA}
I_k & := & \sum_{k_1\ge k} \sum_{|m|\lesssim 2^{k_1-k}} 
\sum_{l\ge 0}   2^k \Bigl| \Bigl(   {\mathcal F}_{xt} \Bigl[A_{3}(\Delta_{k_1} v_1^{k_1,m},v_2^{k_1,m},v_3^{k_1,m})\Bigr], \nonumber  \\
 & & \hspace*{30mm}  \langle \tau-\xi^3+i 2^k \rangle^{-1}    \eta_l (\tau-\xi^3)
\tilde{\eta}_k   \widehat{w}\Bigr)_{L^2}
   \Bigr| \nonumber \\
   & \lesssim&  2^{-\frac{k}{16}}  \sup_{l} (2^{-l/2} \| \nu_j \tilde{\eta}_k  \widehat{w}\|_{L^2})  \| 
   \sup_{m\in \Z}\prod_{i=1}^3  \|S_{k_1+1}  v_i^{k_1,m} \|_{X^{0,1/2,1}_\lambda}  \; .\label{dual3}
  \end{eqnarray}
  Indeed, this will be acceptable, since by proceeding as in \eqref{rg1}-\eqref{rg2}, it is easily checked that 
  \begin{eqnarray*}
  \|S_{k_1+1}  v_i^{k_1,m} \|_{X^{0,1/2,1}_\lambda}^2 
&   \le & \sum_{j=0}^{k_1-1} \|\Delta_j  u_i \|_{F^{1/2}_{\lambda,j,m 2^{-k_1}}}^2+\sum_{q=-1}^1 \|\Delta_{k_1} Êu_i \|_{F^{1/2}_{\lambda,k_1,m 2^{-k_1}+q 2^{-(k_1+1)}}}^2
\\
& & +\sum_{q=-2}^2 \|\Delta_{k_1+1} Êu_i \|_{F^{1/2}_{\lambda,k_1+1,m 2^{-k_1}+q 2^{-(k_1+2)}}}^2
  \lesssim    \|u_i \|_{G_\lambda}^2 \; .
  \end{eqnarray*}
   Since the norms in the right-hand side of \eqref{dual3}Ê only see the size of the modulus of the 
          Fourier transform of the functions we can assume that all the functions have non negative Fourier transforms. 
 We will  use the following notations :
    \begin{equation}\label{defnu}
     {\tilde \nu}_l(\tau,\xi) :=\eta_l(\tau-\xi^3) \mbox{ for } l> k_1 \; \mbox{ and } \quad  {\tilde \nu}_l(\tau,\xi) :=\eta_{\le k_1}(\tau-\xi^3) \mbox{  for } l=k_1\; . 
     \end{equation}
     In view of \eqref{31}, in this region  it holds $ \Bigl| |\xi_1|-|\xi_3|\Bigr|Ê\gtrsim |\xi_1| $ and  $ |\xi-\xi_2|\gtrsim |\xi_1| $. 
  Lemma \ref{lem2}, \eqref{32} and  \eqref{gh1}-\eqref{gh2} thus lead to 
  \begin{eqnarray}\label{est333}
I_k & \lesssim  &  \sum_{k_1\ge k} \sum_{|m|\lesssim 2^{k_1-k}} 
\sum_{l\ge 0}    2^k \Bigl\|Ê \Lambda[2^{k_1}]\Bigl(\eta_{ k_1} \widehat{v_1^{k_1,m}} \,,\,  \eta_{\le k_1+1} \widehat{v_3^{k_1,m}}\Bigr)
\Bigr\|_{L^2}Ê\nonumber \\
 & & \hspace*{20mm} \| \Bigl(\langle \tau-\xi^3+i 2^k \rangle^{-1} \eta_{k} \widehat{w} \Bigr) \ast  \eta_{\le k_1+1} \check{\widehat{v_2^{k_1,m}}}\|_{L^2(|\xi|\gtrsim 2^{k_1})} \nonumber \\
 & &   \lesssim \sum_{k_1\ge k} \sum_{|m|\lesssim 2^{k_1-k}} \sum_{\min(l_1,l_2,l_3) \ge k_1, l\ge 0 \atop \max(2^l,2^{l_i})\ge \lambda^{-1} 2^{2k_1}}   Ê 2^k 2^{(l_1\wedge l_3)/2}  (2^{(l_1\vee l_3)/2}2^{-k_1}  +\lambda^{-1/2})^{3/4} Ê\nonumber \\
 & &  (2^{(l_1\vee l_3)/4}2^{-k_1/4}  +\lambda^{-1/2})^{1/4} 2^{(l\wedge l_2)/2}  (2^{(l\vee l_2)/4}2^{-k_1/4}  +\lambda^{-1/2})   2^{-(l\vee k)} 
 \nonumber \\
   & &  \| \tilde{\nu}_{l_2} \eta_{\le k+1}   \widehat{v_{2}^{k_1,m}}\|_{L^2}  \|\tilde{\nu}_{l_3}  {\eta}_{\le k_1+1}  \widehat{v_{3}^{k_1,m}}\|_{L^2} Ê 
    \|  \tilde{\nu}_{l_1}  \eta_{k_1}    \widehat{v_{1}^{k_1,m}}\|_{L^2}  \| \nu_{l}  \eta_k  \widehat{w}\|_{L^2}
  \nonumber \\
   & \lesssim   &  \sum_{k_1\ge k} \sum_{|m|\lesssim 2^{k_1-k}} 2^k 2^{-\frac{17}{16} k_1}  
     \sup_{l}  \Bigl(2^{-\frac{17}{32}l}   \| \nu_{l}  \eta_k  \widehat{w}\|_{L^2}\Bigr)
     \prod_{i=1}^3 \sup_{l\ge k_1}  \Bigl(2^{\frac{15}{32} l}   \| \tilde{\nu}_{l}  \eta_{\le k_1+1}   \widehat{v_{i}^{k_1,m}}\|_{L^2} 
\Bigr)\nonumber\\
& \lesssim   &  \sum_{k_1\ge k}  2^{-\frac{k_1}{16}} 
     \sup_{l}  \Bigl(2^{-\frac{17}{32}l}   \| \nu_{l}  \eta_k  \widehat{w}\|_{L^2}\Bigr) 
     \sup_{|m|\lesssim 2^{k_1-k}}\Bigl[ 
     \prod_{i=1}^3\|S_{k_1+1} v_{i}^{k_1,m}\|_{X^{0,15/32,1}_\lambda}\Bigr] \label{zq}
         \end{eqnarray}
         which yields \eqref{dual3} by summing in $ k_1 $  and concludes the proof of the lemma. \qed
\begin{corollary}\label{coro1}
 There exists $ \varepsilon_0>0 $ such that  any  $ \lambda\ge 1 $ and any solution $ u \in C^\infty(\R;H^\infty(\lambda \T)) $ to \eqref{mKdV*} satisfying $ \|Êu_0\|_{L^2_\lambda}Ê\le \varepsilon_0   $,  it holds
\begin{equation}\label{estL4} 
\| u\|_{G_\lambda} \lesssim \|Êu_0\|_{L^2_\lambda} \; .
\end{equation}
\end{corollary}
\proof   We are going  to implement a continuity argument on the space period. Recall  that if $ u(t,x) $ is a smooth global  $2\lambda \pi $-periodic solution of \eqref{mKdV*} with initial data
 $ u_0 $ then $ u_\beta(t,x)=\beta^{-1} u(\beta^{-3}t,\beta^{-1} x)$ is a $ (2\pi \lambda\beta) $-periodic solution of 
  \eqref{mKdV*} 
 emanating from $u_{0,\beta}=\beta^{-1} u_0(\beta^{-1} x) $.  Moreover, 
 $$
 \|u_{0,\beta}\|_{L^2_{\lambda\beta}}= \beta^{-1/2} \|u_0 \|_{L^2_\lambda} \mbox{ and  }
  \|u_{0,\beta}\|_{H^1_{\lambda\beta}}\le  \beta^{-1/2} \|u_0\|_{H^1_\lambda}\; .
  $$
From the conservation of the $ L^2 $-norm and of the Energy,
$$E(u):= \frac{1}{2} \int_{\T}u_x^2 \mp \frac{1}{6} \int_{\T} u^4 \; ,
$$
 and Sobolev inequalities (in the focusing case) , 
 we get 
$$\|u_\beta \|_{L^\infty(\R; L^2_{\lambda\beta})}= \|u_{0,\beta}  \|_{L^2_{\lambda\beta}}   \mbox{ and  } 
\|u_\beta \|_{L^\infty(\R; H^1_{\lambda\beta})}\lesssim  \|u_{0,\beta}  \|_{H^1_{\lambda\beta}} (1+\|u_{0,\beta}\|_{L^2_{\lambda \beta}}^2 )
$$
In particular, it follows from \eqref{estF} that for $ \beta \ge 1 $, 
\begin{equation}\label{estF8}
\|u_\beta\|_{G_{\lambda\beta}} \le \, C  \, \Bigl(  \|u_0\|_{L^2_\lambda} +Ê \|u_{\beta}\|_{G_{\lambda\beta}}^3\Bigr) \; ,
\end{equation}
for some constant $ C>0 $.
Now,  from classical linear estimates in Bourgain's spaces and the Duhamel formula,  it holds 
\begin{eqnarray*}
\|u_\beta\|_{G_{\lambda\beta}}  & \lesssim & \sup_{t_0\in\R} \| \eta_0(t-t_0) u_\beta \|_{X^{0,1/2,1}_{\lambda\beta}} 
\lesssim \|Êu_\beta(t_0)\|_{L^2_{\lambda\beta}}Ê+  \|Ê(u_\beta^2-P_0(u_\beta^2))\partial_x u_\beta \|_{L^{2}(]t_0-2,t_0+2[; L^2_{\lambda\beta})}   \\
& \lesssim  & \|Êu_\beta\|_{L^\infty(\R; L^2_{\lambda\beta})}Ê+Ê\|u_\beta\|^3_{L^\infty(\R; H^1_{\lambda\beta})} 
\lesssim \beta^{-1/2} \|u_0 \|_{L^2_\lambda}+\beta^{-3/2} \|u_0\|_{H^1_\lambda}^3 \; .
\end{eqnarray*}
Therefore for $ \beta \ge 1 $ large enough (depending on $\|u_0\|_{H^1}$), $ \|u_\beta\|_{G_{\lambda\beta}} \le  C\, \|Êu_0\|_{L^2_\lambda}$. Recalling that thanks to the infinite number of conservation laws, $ u $ actually belongs to 
 $C^\infty_b(\R;H^\infty(\lambda\T)) $,   $ \beta\mapsto  \|u_\beta\|_{G_{\lambda\beta}} $ is  continuous on $ \R_+^* $ and 
the result follows from   \eqref{estF8} and a  classical continuity argument. 
\section{Proof of Theorem \ref{Existence mKdV}}
In this section we follow the process proposed in  \cite{M1} to identify the limit of a $ L^2$-bounded sequence of solutions to mKdV. 
To pass to the limit on the nonlinear term in \eqref{mKdV*}, we will make use of the space $ F_{\lambda,T}^{s,b} $ introduced in \cite{IKT}. We will need the following lemma which states that for any $ s<0 $, $ b<1/2 $ and $ T>0 $, 
  $ G_\lambda $ is compactly embedded in  $ F^{s,b}_{\lambda,T} $. 
\begin{lem}\label{equicontinuous}
Let $ \lambda \ge 1 $ and $ \{u_n\}_{n\in\N} $ be a bounded sequence of $ G_\lambda $. Then 
 for any $T>0 $, $ s<0 $ and $b<1/2 $, $ \{u_n\}_{n\in \N} $ is relatively compact in $ F^{s,b}_{\lambda,T} $.
\end{lem}
\proof 
First we observe that for any $ s<0 $ and any $ u\in G_\lambda $,
\begin{eqnarray}
\|u\|_{F_\lambda^{s,1/2}}^2 & = & \sum_{k\in \N}Ê2^{2ks} \sup_{t\in \R} \| \Delta_k u \|^2_{F^{1/2}_{\lambda,k,t}} \nonumber \\
& \lesssim & (\sum_{k\in \N} 2^{2ks}) \sup_{k\in \N}Ê\sup_{t\in \R} \| \Delta_k u \|^2_{F^{1/2}_{\lambda,k,t}}\nonumber \\
 & \le & \Bigl( C(s) \| u\|_{G_\lambda}Ê \Bigr)^2 \; .\label{qq1} 
\end{eqnarray}
Hence $ G_\lambda \hookrightarrow  F^{s,1/2}_\lambda $. Second, proceeding as in Remark \ref{remark2}, it is easy to check that there exists $ C_0>0 $ such that for all $ s<0 $ and $ v\in F^{s,1/2}_{\lambda} $, 
\begin{equation}\label{qq2}
\| v\|_{L^\infty(\R;H^s_\lambda)} \le C_0 \, \|Êv\|_{F_\lambda^{s,1/2}} \; .
\end{equation}
Let us now prove the desired result. We will use a diagonal extraction argument. Let us fix $ s'<0 $, $ b'<1/2$ and $ T>0 $. We notice that for any fixed $ k_0\in \N $, 
 there exists $ C_{k_0} >0 $ such that 
 $$
 \| S_{k_0} v \|_{F_{\lambda,T}^{s',b'}}\lesssim \|ÊS_{k_0} v \|_{X^{s',b'}_{\lambda,T}} \le C_{k_0}Ê \| S_{k_0} v \|_{F_{\lambda,T}^{s',b'}}
 $$
 Indeed, the first inequality is obvious and the second one follows from the following chain of inequalities for any $ v\in F^{s',b'}_{\lambda,T}$,
 \begin{eqnarray*}
  \|ÊS_{k_0} v \|_{X^{s',b'}_{\lambda,T}}^2 & \le & \sum_{k=0}^{k_0+1} 2^{2ks'} \|\eta_0(\frac{t}{2T}) \Delta_k S_{k_0} {\tilde v} \|^2_{X^{0,b'}_{\lambda}}\\
   &\lesssim &  \sum_{k=0}^{k_0+1} 2^{2k} T^2 \sup_{t\in \R} 2^{2ks'} \| \Delta_k S_{k_0} {\tilde v} \|^2_{F^{b'}_{\lambda,k,t}}\\
   & \lesssim & \Bigl(C_{k_0}\| S_{k_0} v \|Ê_{F^{s',b'}_{\lambda,T}}\Bigr)^2
 \end{eqnarray*}
 where $ {\tilde v} $ is an extension of $ v $ such that $\| S_{k_0} v \|Ê_{F^{s',b'}_{\lambda}} \le 2  \| S_{k_0} v \|Ê_{F^{s',b'}_{\lambda,T}} $.  Since, according to \eqref{qq1}, $ \{ u_n\}_{n\ge 0}Ê$ is bounded in $ F^{s',1/2}_{\lambda} $, it follows that $ \{ S_{k_0} u_n\}_{n\ge 0} $ 
   is bounded in $ X^{s',b'}_{\lambda} $ for any $ k_0\ge 0 $. Now, using that for $ s<s'$ and $b<b' $, $ X^{s',b'}_{\lambda,T+1} $ is compactly embedded in 
    $ X^{s,b}_{\lambda,T+1} $, we deduce that there exists a subsequence $ \{u_{n_q}\} $ of $ \{ u_n\} $ and a sequence $ \{w_k\} \subset F^{s,b}_{\lambda,T+1}Ê$ such that for any $ k\in \N $,
    \begin{equation}\label{qq4}
    S_k  u_{n_q} \to w_k \mbox{ in } F^{s,b}_{\lambda,T+1} \; .
    \end{equation}
    We define $ w\in {\mathcal S}'(]-T-1,T+1[\times \T) $ by $ \Delta_k w =\Delta_k w_{k+1}Ê$ for all $ k\in \N $. Clearly, for all $ k_0\ge 0 $, 
  $$
   \sum_{k=0}^{k_0} 2^{2ks} \sup_{t\in -]T,T[} \| \Delta_k w \|^2_{F^{b}_{\lambda,k,t}} \lesssim \sum_{k=0}^{\infty}Ê2^{2ks} \sup_{n\in \N}\sup_{t\in \R}Ê\| \Delta_k u_n\|^2_{F^{b}_{\lambda,k,t}}
   \lesssim \sup_{n\in\N}  \|u_n\|_{G_\lambda}^2\; .
  $$
  which ensures that $ w\in F^{s,b}_{\lambda,T} $. Moreover, proceeding as in \eqref{qq1}, it easy to check that 
  $$
  \lim_{k_0\to \infty} \sup_{n\in\N} \| \sum_{k=k_0}^\infty \Delta_k u_n \|_{F^{s,b}_{\lambda,T}}=0 \; .
  $$
 It thus  follows from \eqref{qq4} that 
   $ \| u_{n_q} -w \|_{F^{s,b}_{\lambda,T}} \to 0 $ as $ q\to \infty $. \vspace{2mm} \qed\vspace{2mm} \\
  The following proposition states that our sequence of solution $ \{u_n\} $ is uniformly equi-continuous in $ H^s(\T) $ for any $ s<0 $.  \begin{prop}\label{prolo51}
  Let $ \lambda \ge 1 $ and $ \{u_n\}_{n\in\N} $ be a sequence of smooth solutions to \eqref{mKdV*} that is bounded in  $ G_\lambda $. Then, for  any $ \varepsilon>0 $ and $ s<0 $, there exists $ \delta_{\varepsilon,s} >0 $ such that $ \forall (t_1,t_2)\in \R^2 $, 
$$
|t_1-t_2| < \delta_{\varepsilon,s} \Rightarrow \| u_n(t_1)-u_n(t_2) \|_{H^s_\lambda} < \varepsilon , \; \forall n\in \N \, .
$$
 \end{prop}
  \proof 
  First we claim that if we prove that $\{u_n\} $ is bounded in $ F^{s,b}_\lambda $ for some $b>1/2$ and $ s<0 $ then we are done. 
  To prove this claim,  we fix $ s:=s_0<0 $ and we set 
  $$  M:= \sup_{n\in \N} (\| u_n\|_{G_\lambda}+\|u_n\|_{F^{s_0,b}_\lambda} )\; .
  $$
   For any given $ \varepsilon>0 $, we denote by $ k_{\varepsilon} $ the smaller integer such that 
\begin{equation}\label{qq3} 
\Bigl( \sum_{k=k_{\varepsilon}}^{\infty} 2^{2ks_0}Ê\Bigr)^{1/2} < \frac{\varepsilon}{4C_0\, C(s_0) M } \; .
\end{equation}
From \eqref{qq1}-\eqref{qq2}  and \eqref{qq3} we infer that 
$$
\Bigl\| \sum_{k=k_{\varepsilon}+1}^{\infty}\Delta_k u_n(t_1) - \sum_{k=k_{\varepsilon}+1}^{\infty}\Delta_k u_n(t_2)\Bigr\|_{H^{s_0}_\lambda}Ê
 \le 2 C_0 \, \Bigl\| \sum_{k=k_{\varepsilon}+1}^{\infty}\Delta_k u_n\Bigr\|_{F^{s_0,1/2}_\lambda}Ê\le \varepsilon/2 \, .
 $$
Now, let $\gamma \, :\,  \R\to [0,1] $ be a $ C^\infty$-function  with compact support in $[-1,1] $ satisfying $ \gamma\equiv 1 $ on 
 $[-1/2,1/2] $. By \eqref{ds3}, for $ |t_1-t_2|< 2^{-k_{\varepsilon}} /4 $ and any $ R>0 $, it holds 
\begin{eqnarray}\arraycolsep2pt
\Bigl\| \sum_{k=0}^{k_{\varepsilon}} \Bigl( \Delta_k u_n(t_1)&  - & \Delta_k u_n(t_2)\Bigr)  \Bigr\|_{H^{s_0}_\lambda}^2\nonumber \\
 & =&  
 \Bigl\|\sum_{k=0}^{k_{\varepsilon}}\Bigl(  \gamma(2^k (t_1-t_1))  \Delta_k u_n(t_1) - \gamma(2^k( t_2-t_1)) \Delta_k u_n(t_2)\Bigr)  \Bigr\|_{H^{s_0}_\lambda}^2\nonumber\\
  & \sim  &  \sum_{k=0}^{k_{\varepsilon,s}}  2^{2ks_0}\Bigl\|    \int_{\R} {\mathcal F}_{t,x} \Bigl(  \gamma(2^k(t-t_1))\Delta_k u_n\Bigr)(\xi,\tau) 
  ( e^{it_1\tau}-e^{it_2 \tau} )\, d\tau    \Bigr\|_{L^2_\xi}^2\nonumber\\ 
  & \lesssim&  C(s_0) \sup_{0\le k \le k_{\varepsilon} }\Bigl[ R^2  |Êt_1-t_2|^2Ê \Bigl\| \int_{-R}^RÊ\Bigl|   {\mathcal F}_{t,x} \Bigl( \Delta_k \gamma(2^k(t-t_1)) u_n\Bigr)(\xi,\tau) \Bigr|
  \, d\tau \Bigr\|_{L^2_\xi}^2 \nonumber \\
  & & +  \Bigl\| \int_{|\tau|>R}Ê\Bigl|   {\mathcal F}_{t,x} \Bigl( \Delta_k \gamma(2^k(t-t_1)) u_n\Bigr)(\xi,\tau) \Bigr|
  \, d\tau \Bigr\|_{L^2_\xi}^2 \Bigr] \nonumber \\
   & \lesssim & C(s_0)  R^2  |Êt_1-t_2|^2 \sup_{0\le k\le k_{\varepsilon}} \Bigl[ \Bigl\| \Delta_k \gamma(2^k(t-t_1)) u_n \Bigr\|_{X^{0,1/2,1}_\lambda}^2 \nonumber\\
    & &+ R^{1-2b}\Bigl\| \Delta_k \gamma(2^k(t-t_1)) u_n \Bigr\|_{X^{0,b,1}_\lambda}^2\Bigr]\nonumberÊ\\
    & \lesssim& C(s_0) M\Bigl( |t_1-t_2|^2 R^2+ R^{1-2b}2^{(2b-1)k_{\varepsilon}}\Bigr)  \; ,\label{ssq}
\end{eqnarray}
\arraycolsep4pt
where in the last step we use that, according to \eqref{ds1}Ê-\eqref{ds2} , for all $ b\in [1/2,1[ $ and all  $k\in \N  $, 
$$
 \|\gamma(2^k t) f\|_{X^{0,b,1}_\lambda} \lesssim  2^{(b-\frac{1}{2}) k}Ê\|Êf\|_{X^{0,b,1}_\lambda} \, , \quad \forall f \in X^{0,b,1}_\lambda \; .
 $$
 Taking  $ R:= |t_1-t_2|^{-\frac{1}{b+1/2}} $, \eqref{ssq} leads to 
 $$
 \Bigl\| \sum_{k=0}^{k_{\varepsilon}} \Bigl( \Delta_k u_n(t_1)- \Delta_k u_n(t_2)\Bigr)  \Bigr\|_{H^{s_0}_\lambda}
 \lesssim C(\varepsilon,M,b)  |t_1-t_2|^{\frac{b-1/2}{b+1/2}}
 $$
This  gives the desired result for $ s=s_0 $ and   $ \delta_{\varepsilon,s_0}:=
 \min [\frac{\varepsilon}{2M}, (\frac{\varepsilon}{2C(\varepsilon,M,b)})^\frac{b+1/2}{b-1/2}]$.  Finally, the result for any $ s\in ]s_0,0[Ê$ follows by interpolating with \eqref{lll}.\vspace{2mm}
 
 We are thus reduced  to prove that $\{u_n\} $ is bounded in $ F^{-2,b}_{\lambda} $ for some $ b>1/2$.  We proceed as in Proposition \ref{prolo} and write that 
    for any $( t_0,t)  \in \R^2 $, it holds 
 $$
  u(t)= U(t-t_0)  u(t_0)+\int_{t_0}^t U(t-t') \partial_x \Bigl(A({u}(t'))+B({u}(t')\Bigr)dt' \; .
 $$
 By translation in time we can always assume that $ t_0=0 $ and 
 according to Lemmas \ref{linear1}-\ref{linear2},  for any $ b\in ]1/2,1[ $, 
 $$
 \|U(t)\Delta_k u(0)\|Ê_{F^{b}_{\lambda,k,0}}Ê\lesssim 2^{(b-1/2)k} \| \Delta_k u (0)\|_{L^2_\lambda}
 Ê\lesssim 2^{(b-1/2)k} \| \Delta_k u\|_{L^\infty_t L^2_\lambda}
 $$
 and 
 $$
 \Bigl\|\int_{0}^t U(t-t') \partial_x \Delta_k (A(u(t'))+B(u(t'))dt' \Bigr\|_{F^{b}_{\lambda,k,0}} \lesssim
2^{(b-\frac{1}{2})k} \, 2^k   \| \Delta_k (A(u))+B(u)) \|_{Z^{b}_{\lambda,k,0}}\; .
 $$
Therefore, for any $ b\in ]1/2,1[ $,
\begin{eqnarray}
\|U(t) u\|_{F^{-2,b}_\lambda}^2 &=&\sum_{k\in \N} 2^{-4k} \sup_{t\in\R} \|U(t)\Delta_k u \|Ê_{F^{b}_{\lambda,k,t}}^2
\nonumber \\
& \lesssim  & \sum_{k\in \N} 2^{-4k} 2^{(2b-1)k} \| \Delta_k u\|_{L^\infty_t L^2_\lambda}^2\nonumber \\
& \lesssim & C(b)  \| u\|_{L^\infty_t L^2_\lambda}^2\lesssim C(b) \|u\|_{G_\lambda}^2 \; ,
\end{eqnarray}
where we used \eqref{lll} in the last step. In the remaining we assume that $ b\in ]1/2,1/2+\frac{1}{32}[ $.  It remains to prove that 
 $$
 \sum_{k=0}^\infty  \sup_{t\in \R} \Bigl(  \| \Delta_k  A({u})\|_{Z^b_{\lambda, k,t}}^2+
 \|  \Delta_k B(u) \|_{Z^b_{\lambda, k,t}}^2\Bigr)  \lesssim  \| u\|_{G_\lambda}^6\; .
 $$
To control the term involving $ B $, we proceed as in Lemma \ref{bil2}. By duality it suffices to prove that 
\begin{eqnarray*}
{\tilde I}_k& := &  \Bigl| \Bigl(  {\mathcal F}_{xt} \Bigl[ B(\Delta_k v_1,{\tilde \Delta}_k v_2,{\tilde \Delta}_k v_3) \Bigr],
 \langle \tau-\xi^3+i 2^k \rangle^{-1}     \widehat{w} \Bigr)_{L^2}
   \Bigr| \nonumber \\
   & \lesssim& 2^{-k/2}   \sup_{l} (2^{-bl} \| \nu_l \widehat{w}\|_{L^2})  
     \prod_{i=1}^3 \|S_{k+1} v_i\|_{X^{0,1/2,1}_\lambda}  \; .
  \end{eqnarray*} 
where $ v_i:=\eta_0(2^k t) u_i $. For $ k=0,1,2, $ this follows directly from \eqref{strichartz}, whereas for $ k\ge 3 $ we deduce from \eqref{kkk}  that 
$$
{\tilde I}_k \lesssim 2^{-k} 2^{(b-\frac{1}{2})k}Ê\sup_{l} (2^{-bl} \| \nu_l \widehat{w}\|_{L^2})  
     \prod_{i=1}^3 \|S_{k+1} v_i\|_{X^{0,1/2,1}_\lambda} 
     $$
     which is acceptable.
     
      Now to estimate the term involving $ A$ we denote by $ \xi_i  $ the Fourier modes of $ u_i $ and we assume  by symmetry that $ |\xi_1|\ge |\xi_2|\ge |\xi_3 |$. 
 We divide $ A $ in different terms corresponding to regions of  $ (\lambda \Z)^3 $. \\
 {\bf 1.}Ê{$ |\xi_1|\le 2^4$}. This region can be treated as in Lemma \ref{estimA}Ê by using \eqref{strichartz}.\\
 {\bf 2.}Ê$|\xi_1|\ge 2^4$ and   $ |\xi_1|\le 4  |\xi|$. By duality it suffices to prove that 
  \begin{eqnarray*}
 {\tilde H}_{k}& := & 2^k \Bigl| \Bigl( {\mathcal F}_{xt} \Bigl[ A_1(v_1, v_2, v_3) \Bigr],
 \langle \tau-\xi^3+i 2^k \rangle^{-1}     \widehat{\tilde{\Delta} w} \Bigr)_{L^2}
   \Bigr| \nonumber \\
   & \lesssim&  2^{-k/2}  \sup_{l} (2^{-bl} \| \nu_l  \widehat{w}\|_{L^2})   
     \prod_{i=1}^3   \|S_{k+3} v_i\|_{X^{0,1/2,1}_\lambda} 
       \end{eqnarray*} 
where $ A_1 $ is defined in \eqref{defA1} and $ v_i:=\eta_0(2^k t) u_i $. 
This can be done by proceeding  as above for $ B $, separating as in Lemma \ref{estimA}, the 4 cases :
 $\xi\xi_1\le 0$, $ \xi \xi_2\le 0 $, $\xi\xi_3\le 0 $ and all the $ \xi_i $'s  of the same sign. \\
 {\bf 3.}Ê$|\xi_1|\ge 2^4 $ and $ |\xi_1|\ge 4 |\xi|$. In this last  region it suffices to prove that 
\begin{eqnarray*}
{\tilde J}_{k}  &:=& Ê \sum_{k_1\ge k} 
\sum_{l\ge 0} 2^{bl}  2^k \Bigl\| \eta_l (\tau-\xi^3)\langle \tau-\xi^3+i 2^k \rangle^{-1}  
 \tilde{\eta}_k(\xi){\mathcal F}_{tx} \Bigl(A_{3}(\Delta_{k_1} v_1,v_2,v_3)\Bigr) \Bigr\|_{L^2} \\
 & \lesssim & 2^{-\frac{k}{2}} \prod_{i=1}^3 \|u_i\|_{G_\lambda}
\end{eqnarray*}
where $ A_3 $ is defined in \eqref{defA3}. But this follows directly from \eqref{zq} for any $ b\in ]\frac{1}{2}, \frac{1}{2}+\frac{1}{32}[ $. This completes the proof of the lemma. \qed\vspace{2mm} 

As noticed in \cite{M1}, $ B $ has got a nice structure for passing to the limit in the sense of distributions. More precisely, we have the following lemma:
\begin{lem}\label{lemB}
For any $ \lambda\ge 1 $ and $ T>0 $, the operator $ \partial_x B $ is continuous from $(F^{-1/3,1/3}_{\lambda,T})^3 $ into $ X^{-4,-1/3}_{\lambda,T} $.
\end{lem}
\proof 
 Taking $ w\in  X^{4,1/3}_{\lambda} $ supported in time in $ ]-2T,2T[  $ and extensions $ {\tilde u_i}\in
F^{-1/3,1/3}_\lambda $ of $ u_i \in F^{-1/3,1/3}_{\lambda,T} $   such that $
\|{\tilde u_i}\|_{F^{-1/3,1/3}_\lambda}
 \le 2 \|u_i\|_{F^{-1/3,1/3}_{\lambda,T}} $,
 it holds
\begin{eqnarray*}
I& := & \Bigl| \Bigr(w, \partial_x B(u_1,u_2,u_3)\Bigl)_{L^2(\R\times\lambda \T)} \Bigr| \\
& \lesssim & \sum_{k=0}^\infty \Bigl|  \Bigr(\Delta_k w_x , B({\tilde \Delta}_k {\tilde u_1},{\tilde \Delta}_k {\tilde u_2}, {\tilde \Delta}_k {\tilde u_3})\Bigl)_{L^2(\R\times\lambda \T)} \Bigr| \\
& \lesssim &  \sum_{k=0}^\infty \sum_{|m|\lesssim T 2^k}Ê\Bigl|  \Bigr(\Delta_k w_x , B({\tilde \Delta}_k v_1^{k,m_1},{\tilde \Delta}_k v_2^{k,m_2}, {\tilde \Delta}_k v_3^{k,m_3})\Bigl)_{L^2(\R\times\lambda \T)} \Bigr| 
\end{eqnarray*}
 where for any $ k\in \N $, $ i\in \{1,2,3\} $ and $ m \in \Z $, we set $ v_i^{k,m}:=\gamma(2^k t-m) {\tilde u_i} $  with $ \gamma $ defined as in \eqref{defgamma}. 
 \eqref{strichartz} and Bernstein inequality then ensure that 
 \begin{eqnarray}
 I & \lesssim & \sum_{k=0}^\infty \sum_{|m|\lesssim T 2^k} 2^{2k} \| {\mathcal F}^{-1} (|\widehat{\Delta_k w} |) \|_{L^4_\lambda} 
 \prod_{i=1}^3 2^{-k/3} \| {\mathcal F}^{-1} (|\widehat{{\tilde \Delta_k} v_i^{k,m}} |) \|_{L^4_\lambda} \nonumber\\ 
 & \lesssim & T  \sum_{k=0}^\infty 2^{3k} \| \Delta_k w \|_{X^{0,1/3}_\lambda} \prod_{i=1}^3 \sup_{m\in \Z} 
  \| {\tilde \Delta}_k  v_i^{k,m}\|_{X^{-1/3,1/3}_\lambda}\label{uuy}
 \end{eqnarray}
 But on account of \eqref{gh1} (with obvious modification for $ k=0$), 
 \begin{eqnarray*}
 \| {\tilde \Delta}_k v_i^{k,m} \|_{X^{-1/3,1/3}_\lambda} & \lesssim & 
  \sum_{j=-1}^1 \| \Delta_{k+j} \gamma(2^k t-m) {\tilde u_i} \|_{X^{-1/3,1/3}_\lambda}\\
  & \lesssim & \sum_{j=-1}^1 2^{-(k+j)/3} \sup_{t\in \R} \|  \Delta_{k+j}{\tilde u_i}Ê\|_{F^{1/3}_{\lambda,k+j,t}}\\
&   \lesssim & \| {\tilde u_i}\|_{F^{-1/3,1/3}_\lambda}\lesssim \| u_i\|_{F^{-1/3,1/3}_{\lambda,T}}Ê\; .
 \end{eqnarray*}
 Therefore, \eqref{uuy} leads to 
 \begin{equation}\label{momo1}
 I \lesssim T \| w\|_{X^{4,1/3}_\lambda} \prod_{i=1}^3 \| u_i\|_{F^{-1/3,1/3}_{\lambda,T}}
 \end{equation}
 which concludes the proof of the lemma. \qed
 
Let us now prove a continuity result for the non resonant part $ A$.
\begin{lem} \label{lemA}
For any $ \lambda\ge 1 $ and $ T>0 $, the operator $  \partial_x A $   is continuous from $(F^{-2^{-6},15/32}_{\lambda,T})^3 $ into $ X^{-4,-1/2}_{\lambda,T} $.
\end{lem} 
\proof 
 Taking $ w\in  X^{4,1/2}_{\lambda}$ supported in time in $]-2T,2T[ $ and extensions $ {\tilde u_i}\in
F^{-2^{-6},15/32}_\lambda $ of $ u_i \in F^{-2^{-6},15/32}_{\lambda,T} $   such that $
\|{\tilde u_i}\|_{F^{-2^{-6},15/32}_\lambda}
 \le 2 \|u_i\|_{F^{-2^{-6},15/32}_{\lambda,T}} $,
 it holds
\begin{eqnarray*}
J& := & \Bigl| \Bigr(w, \partial_x A(u_1,u_2,u_3)\Bigl)_{L^2(\R\times\lambda \T)} \Bigr| \\
& \lesssim & \sum_{(k,k_1,k_2,k_3)\in \N^4 } \Bigl|  \Bigr(\Delta_k w_x , A({ \Delta}_{k_1} {\tilde u_1},{ \Delta}_{k_2} {\tilde u_2}, {\Delta}_{k_3} {\tilde u_3})\Bigl)_{L^2(\R\times\lambda \T)} \Bigr| \\
\end{eqnarray*}
By symmetry we may assume that $ k_1\ge k_2\ge k_3 $. Now, the sum in the reion $ k_1\le 6 k+10 $ can clearly be treated exactly in the same way as 
 we treat $ B $ in the preceding lemma. It thus remains to consider the sum over the region $ k_1>6 k+10 $. We follow  the proof of Lemma \ref{estimA}. First  we notice that \eqref{31} -\eqref{32} Êhold in this region. Then setting,    $$ v_i^{k_1,m}: = \gamma(2^{k_1}t -m) {\tilde u_i} , \;  i=1,2,3 , \, m\in \Z , $$ 
  with $ \gamma $ defined as in \eqref{defgamma},  we obtain in this region
$$
J  \lesssim \sum_{k\ge 0}Ê \sum_{k_1\ge k_2\ge k_3 \atop k_1\ge 6k+10 }  \sum_{|m|\lesssim  T 2^{k_1}} 
\Bigl|  \Bigr(\Delta_k w_x , A({ \Delta}_{k_1} { v_1^{k_1,m}},{ \Delta}_{k_2} { v_2^{k_1,m}}, {\Delta}_{k_3} { v_3^{k_1,m}})\Bigl)_{L^2(\R\times\lambda \T)} \Bigr|
$$
Proceeding exactly as in \eqref{zq}, with \eqref{gh1} in hand,  we get
\begin{eqnarray}
J & \lesssim  & \sum_{k\ge 0}Ê \sum_{k_1\ge k_2\ge k_3 \atop k_1\ge 6k+10 } 2^{-\frac{k_1}{16}}\| \Delta_k w\|_{X^{2,1/2}_\lambda}\prod_{i=1}^3 \sup_{m\in \Z} 
  \| \Delta_{k_i}Ê v_i^{k_1,m}\|_{X^{0,15/32}_\lambda}\nonumber \\
  & \lesssim &  T\sum_{k\ge 0}\sum_{ k_1\ge 0} 2^{-\frac{k_1}{16}}k_1^2 \| \Delta_k w\|_{X^{2,1/2}_\lambda}
   \sup_{m\in \Z}(   \| \Delta_{k_1}Ê v_1^{k_1,m}\|_{X^{0,15/32}_\lambda} ) \prod_{i=2}^3  \sup_{0\le k_3\le k_2\le k_1\atop m\in \Z }
  \| \Delta_{k_i}Ê v_i^{k_1,m}\|_{X^{0,15/32}_\lambda}\nonumber\\ 
   & \lesssim  & T \|w\|_{X^{3,1/2}_\lambda} \prod_{i=1}^3 \sup_{t\in \R, k\in \N } 2^{-2^{-6}k} \| \Delta_k {\tilde u_i} \|_{F^{15/32}_{\lambda,k,t}}\nonumber\\
    & \lesssim  & T  \|w\|_{X^{3,1/2}_\lambda} \prod_{i=1}^3 \| u_i \|_{F^{-2^{-6}, 15/32}_{\lambda,T}}\label{momo2}
  \end{eqnarray}
which  completes the proof of the lemma .\vspace{2mm} 
\qed\\
For any $ T>0 $,  let us define the operator $  \Upsilon_T $ which to $ u $ associates
$$
\Upsilon_T(u):={\mathcal F}^{-1}_{tx}\Bigl( \Gamma_T(u)\Bigr) 
$$
where 
\begin{eqnarray}
\Gamma_T(u)(\xi,\tau) &:= & \frac{6i \xi}{\lambda^2} \int_{\tau_1+\tau_2+\tau_3=\tau}  \Bigl[  \widehat{\psi_T u}(\xi,\tau_1)  \widehat{\psi_T u}(\xi,\tau_2)  \widehat{\psi_T u}(-\xi,\tau_3) \nonumber \\
& &\hspace*{-20mm} -\frac{1}{3} \hspace*{-5mm} \sum_{\xi_1+\xi_2+\xi_3=\xi \atop (\xi_1+\xi_2)(\xi_1+\xi_3)(\xi_2+\xi_3) \neq 0} \hspace*{-10mm}  \widehat{\psi_T u}(\xi_1,\tau_1)  \widehat{\psi_T u}(\xi_2,\tau_2)  \widehat{\psi_T u}(\xi_3,\tau_3)  \Bigr]d\tau_1\, d\tau_2\, d\tau_3 \, ,  \label{defLambda}
\end{eqnarray}
with $\psi_T(\cdot):=\psi(\cdot/T)$ for $ \psi $  defined as in Section \ref{section22}. \\
 In view of \eqref{defAB}, for any smooth function $ u \in {\mathcal S}( \lambda\T\times \R) $, $ \Upsilon_T(u) \equiv 6(u^2-P_0(u^2))u_x $ on 
  $ ]-T,T[$. 
From the two above lemmas we infer that $ \Upsilon_T $ can be continuously extended in  $  F^{-2^{-6},15/32}_{\lambda}$ with values in $ 
 X^{-4,-1/2}_{\lambda} $. In particular, for any $ u\in   F^{-2^{-6},15/32}_{\lambda} $ and any $T>0$, this operator 
  defined an element of $ 
 X^{-4,-1/2}_{\lambda} $ with a  $X^{-4,-1/2}_{\lambda} $-norm which is, according to \eqref{momo1} and \eqref{momo2}, of order at most $ O(T) $. This ensures that we can pass to the limit in $ {\mathcal S}' $ on $ \Gamma_T(u) $ as $ T\to \infty $. We  can  thus define the operator
  $ \Gamma $   from  $F^{-2^{-6},15/32}_{\lambda} $ into $ {\mathcal S}'(\lambda\T\times \R) $ by setting 
 $$
 \langle  \Gamma(u),\phi \rangle_{{\mathcal S}',{\mathcal S}},:=\lim_{T\to \infty}\langle  \Gamma_T(u),\phi \rangle_{{\mathcal S}',{\mathcal S}}, \quad \forall \phi 
 \in   {\mathcal S}(\lambda\T\times \R) 
 $$
 Obviously, $ {\mathcal F}^{-1}_{tx}(\Gamma(u))\equiv 
   6(u^2-P_0(u^2))u_x $  on $\lambda\T\times\R $   for any $ u\in {\mathcal S}  (\lambda\T\times \R) $. 
 \begin{definition}\label{def1}
 We will say that a function $ u $ is a weak solution of \eqref{mKdV*} if it satisfies \eqref{mKdV*} in the sense of distributions, when $(u^2-P_0(u^2))u_x $ is interpreted as the inverse Fourier transform of $ \Gamma(u) $. 
 \end{definition}
\begin{prop}\label{propo}
 Let $ \{u_{0,n}\} \subset H^\infty(\T) $ be  such that $ u_{0,n} \rightharpoonup u_0 $ in $ L^2(\T) $. Then 
 there exist a weak solution $ u\in C_{w}(\R; L^2(\T)) \cap \Bigl( \displaystyle \bigcup_{s<0} F^{s,1/2}\Bigr) $ to \eqref{mKdV*}, with $u(0)=u_0 $,  and  a  subsequence of emanating solutions $\{u_{n_k}\} $ to \eqref{mKdV*}  such that for all $T>0 $ and $ \phi\in L^2(\T) $,
\begin{equation}
(u_{n_k}(t),\phi)_{L^2(\T)}  \to  (u(t),\phi)_{L^2(\T)} \mbox{ in } C([-T,T]) \label{ddd}\; .
\end{equation}
 Moreover, $ P_0(u(t))=P_0(u_0)$ for all $ t\in \R $. 
\end{prop}
\proof We proceed as in \cite{M1}. By Banach-Steinhaus theorem, the sequence $ \{u_{0,n}\}Ê$ is bounded in $ L^2(\T) $. We start by assuming that 
 $ \sup_{n\in \N}Ê\|u_{0,n}\|_{L^2} \le \varepsilon_0 $ so that the conclusions of Corollary \ref{coro1} Êhold. 
  From the conservation of the $ L^2$-norm,  the sequence of emanating solutions $ \{u_n\} $ to \eqref{mKdV*} is bounded in $L^\infty(\R;L^2(\T)) $ and thus,  up to a subsequence,  $ \{u_n\} $  converges weakly star in $ L^\infty(\R;L^2(\T))$ to some $ u\in  L^\infty(\R;L^2(\T))$. In particular,
    $ \{\partial _t u_n\} $ and  $ \{\partial_x^3 u_n\} $ converge in the distributional sense to respectively $ u_t $ and $ u_{xxx} $. It remains to pass to the limit on the nonlinear term. By Corollary  \ref{coro1}, $ \{u_n\} $ is bounded in $ G $. 
       From Lemma \ref{equicontinuous} and lemmas  \ref{lemB}-\ref{lemA} it follows that $ \partial_x (A(u_n) + B(u_n)) $ converges to ${\mathcal F}_{tx}^{-1}(\Gamma(u)) $  in the distributional sense on $ ]-T,T[\times \T $. Since this holds for  all $ T>0 $, the convergence holds actually in the distributional sense on $ \R\times\T $.  Therefore, $ u $ is a weak solution to 
     \eqref{mKdV*}  in the sense of Definition \ref{def1}. Note also that, in view of \eqref{qq1}, $ u\in  \displaystyle \cup_{s<0} F^{s,1/2}$. 
      Moreover, according to Proposition  \ref{prolo51}, for any time-independent $ 2\pi $-periodic smooth function  $ \phi $, the family 
       $ \{t\mapsto (u_n(t),\phi)_{L^2_x} \} $ is bounded and 
      uniformly equi-continuous on $ [-T,T] $. Ascoli's theorem then ensures that $ \{t\mapsto (u_{n_k}, \phi)\} $ converges to $t\mapsto (u,\phi)$ in $ C([-T,T]) $. Since $ \{u_n\} $ is bounded in $  L^\infty(\R;L^2(\T)) $, this convergence also holds for any $ \phi \in L^2(\T) $.  This proves \eqref{ddd} and that $ u\in C_{w}(\R;L^2(\T)) $. In  particular, $ u(0)=u_0 $ and for all $ t\in \R$,
  $$
   \int_{\T} u(t)\, dx=\lim_{k\to \infty} \int_{\T} u_{n_k}(t)\, dx= \lim_{k\to \infty}  \int_{\T} u_{n_k}(0)\, dx=
   \int_{\T} u_{0}\, dx\; . 
   $$
   This concludes the proof of the proposition whenever  $ \sup_{n\in \N}Ê\|u_{0,n}\|_{L^2} \le \varepsilon_0 $. Now, if $ \sup_{n\in \N}Ê\|u_{0,n}\|_{L^2} 
   =M> \varepsilon_0 $ we use the dilation symmetry of \eqref{mKdV*}. Recall  that if $ u(t,x) $ is a smooth global  $2 \pi $-periodic solution of \eqref{mKdV*} with initial data
 $ u_0 $ then $ u^\lambda(t,x)=\lambda^{-1} u(\lambda^{-3}t,\lambda^{-1} x)$ is a $ 2\pi \lambda $-periodic solution of 
  \eqref{mKdV*} 
 emanating from $u_{0}^{\lambda}=\lambda^{-1} u_0(\lambda^{-1} x) $.  Setting $ \lambda_M =M^2/\varepsilon^2 $ so that 
 $$
 \sup_{n\in \N} \|u_{0,n}^{\lambda_M}\|_{L^2_{\lambda_M}}Ê\le \varepsilon_0\; ,
 $$
 it follows from above that the conclusions of the proposition hold if one replaces $ \{u_{0,n}\} $, $ \{u_n\} $ and $ u$ by respectively 
 $ \{u_{0,n}^{\lambda_M}\} $, 
  $ \{u_n^{\lambda_M}\} $ and $ u^{\lambda_M} \in C_{w}(\R; L^2(\lambda_0 \T)) \cap \Bigl( \displaystyle \cup_{s<0} F^{s,1/2}_{\lambda_M}\Bigr) $.  This ensures\footnote{It can be easily checked that the dilation symmetry $ u\mapsto \lambda u(\lambda^3 t,\lambda x) $ is an isomorphism from 
   $ F^{s,1/2} $ into $ F^{s,1/2}_{\lambda} $ for any $ \lambda\ge 1 $ and $ s\in \R $.} that these conclusions also hold for $ \{ u_{0,n} \}$, $Ê\{u_n\}$ and $ u$ with $ u(t,x):=\lambda_M u^{\lambda_M} (\lambda_M^{3} t, \lambda_M  x) $ and completes the proof of the proposition. \qed\vspace{2mm}Ê
  
 {\it Proof of  Theorem \ref{Existence mKdV} :}  Again we have to introduce a notion of weak solution for mKdV: 
 \begin{definition}\label{def2}
 We will say that a function $ v \in C_w(\R;L^2(\T)) $ is a weak solution of \eqref{mKdV} with initial data $ v_0 $ if it satisfies the following equation in the sense of distributions,
 \begin{equation}
 v_t+v_{xxx} \mp 6(v^2-P_0(v^2))v_x \pm 6 P_0(v_0^2) v_x=0 \;, \label{mKdV2}
  \end{equation}
  when $(v^2-P_0(v^2))v_x $ is interpreted as the inverse Fourier transform of $ \Gamma(v) $. 
 \end{definition}
 Of course, any smooth solution of the mKdV equation is a weak solution since the $ L^2$-norm is a constant of the motion for smooth solutions.

 Let $ v_0\in L^2(\T) $. It is easy to construct a sequence  
   $ \{v_{0,n}\}\subset H^\infty(\T) $ such that  $\|v_{0,n}\|_{L^2}=\|v_0\|_{L^2} $ for any $ n\ge 0 $ and  $ v_{0,n}\to v_0 $ in $ L^2(\T) $. The emanating solutions $ v_n $ satisfies 
   $$
    v_{n,t}+v_{n,xxx} \pm 6 P_0(v_0^2) v_x \mp 6(v_n^2-P_0(v_n^2))v_x =0 \; .
   $$
  We will consider $ 6 P_0(v_0^2) v_x$ as a part of the linear group of the equation. The $ v_n$ satisfy the same Duhamel formula as the solution of \eqref{mKdV*} where we substitute the linear group of the KdV equation $ U(t) $ by the linear group $ V(t) $ defined by 
 \begin{equation}\label{taz}
\widehat{V(t)\varphi}(\xi)=e^{i q(\xi) t} \,
\hat{\varphi}(\xi) , \quad \xi\in \lambda^{-1}\Z \quad , \; q(\xi)=\xi^3 \mp 6   P_0(v_0^2) \xi \; .
\end{equation}
It is direct to check that all linear estimates remain true when changing the functional spaces in consequence. Also the bilinear estimates in Lemma \ref{lem2}  remain true since the resonance relation remain unchanged (see the proof of this lemma  in the appendix). Therefore, the results of Section \ref{4.1} and the conclusions of 
   Proposition \ref{propo} remain true when substituting the function spaces associated with $ U(\cdot) $ by those associated with $ V(\cdot) $ and \eqref{mKdV*} by \eqref{mKdV2}. We  thus  
   obtain a weak solution $ v\in C_{w}(\R;L^2(\T)) \cap \Bigl( {\displaystyle \cup_{s<0}}{\tilde  F}^{s,1/2}\Bigr) $ to \eqref{mKdV2}    with initial data $ v_0 $, where  ${\tilde  F}^{s,1/2}$ is defined as $F^{s,1/2}$  but for the group $ V(\cdot)$.   Moreover,  by the weak convergence result \eqref{ddd}, $ \|v(t)\|_{L^2} \le \|v_0\|_{L^2} $  for all $ t\in \R $ and thus the weak continuity of $ v$ ensures that $ v(t)\to v(0) $ in $ L^2(\T) $ as $ t\to 0$. This completes the proof of assertion $ i)$. 

Now in the defocusing case, according to \cite{KT2}, the sequence $ \{v_n\} $ of solutions to mKdV emanating from 
 $ \{v_{0,n}\} $ converges in $ C(\R;L^2(\T)) $ to some function $ w$ such that $ \|w(t)\|_{L^2} =\|v_0\|_{L^2}Ê $ for all $ t\in \R $. By the uniqueness of the limit in $ {\mathcal D}'(\R\times \T) $, $ w\equiv v$ on $ \R $ and thus $ w$ 
 is a weak solution of the defocussing  mKdV equation in the sense of Definition \ref{def2}. Actually, using the conservation of the $ L^2 $-norm for $w$, we also obtain that $ w$ satisfy the following equation in the sense of distributions,
 \begin{equation}\label{mKdV3*}
 w_t+w_{xxx} - 6(w^2-P_0(w^2))w_x + 6 P_0(w^2) w_x=0 \;, 
  \end{equation}
  when $(w^2-P_0(w^2))w_x $ is interpreted as the inverse Fourier transform of $ \Gamma(w) $. 
This  concludes the  proof of assertion $ ii)$. 
 \section{Proof of Theorem \ref{main}}
 \subsection{The mKdV equation}
 We will prove that the solution-map $ u_0 \mapsto u $ is not continuous at any $ u_0\in H^\infty(\T)  $  from $ L^2(\T) $ equipped with its weak topology into ${\mathcal  D}'(]0,T[\times \T) $. This obviously  leads to the desired result since $ L^2(\T) $ is compactly embedded in 
 $ H^s(\T) $ for any $ s<0$. Since the sign in front of the nonlinear term will not play any role in the proof, we choose to take the plus sign to simplify the notations.
 
  Let $ u_0\in H^\infty(\T) $ be a non constant function and $ \kappa\neq 0 $ be a real number. We set 
 $$
 u_{0,n}=u_0+\kappa \cos(nx) 
 $$
so that $ u_{0,n} \rightharpoonup u_0 $ in $ L^2(\T) $ and $ \|u_{0,n}\|^2_{L^2} \to  \|u_{0}\|^2_{L^2} +\kappa^2 \pi$. According to Proposition \ref{propo} there exists a subsequence $ \{u_{n_k}\} $ of  the emanating solutions $ \{u_n\} $
 to \eqref{mKdV*}  and $ u\in C_{w}(\R;L^2(\T)) $ a weak solution of  \eqref{mKdV*} , with $ u(0)=u_0 $, satisfying 
$ u_{n_k}(t)   \rightharpoonup u(t) $ in $ L^2(\T) $ for all $ t\in \R $. Let  now $ \{v_{n_k}:=u_{n_k}(\cdot,\cdot-\frac{6t}{2\pi} \|u_{0,n_k}\|_{L^2}^2) \} $ be  the associated subsequence  of solutions to mKdV emanating from $ 
\{u_{0,n_k}\} $. We proceed by contradiction. Assuming  that the solution-map is continuous at $ u_0 $  from $ L^2(\T) $ equipped with its weak topology into ${\mathcal  D}'(]0,T[\times \T) $, we obtain that $ \{v_{n_k} \} $ converges in the sense of distributions in 
$ ]0,T[\times \T $ to the solution  
$v\in C^\infty(\R; H^\infty(\T))$ of mKdV emanating from $ u_0$. It follows  that $\{u_{n_k}\} $ converges in the same sense to 
$ v\Bigl( \cdot, \cdot+\frac{6t}{2\pi}  (\|u_{0}\|_{L^2}^2+\kappa^2 \pi)\Bigr)$ and thus 
\begin{equation}
u\equiv v\Bigl( \cdot, \cdot+\frac{6t}{2\pi}  (\|u_{0}\|_{L^2}^2+\kappa^2 \pi)\Bigr)  \mbox{ on }Ê
]0,T[ \; .  \label{bc}
\end{equation}
This ensures that $ u $ is actually a strong  solution of \eqref{mKdV*} and  satisfies this equation everywhere on 
 $]0,T[\times \T $. On the other hand, according to \eqref{bc}, $ u$ is also solution of 
$$
u_t +u_{xxx} + 6\Bigl(u^2-P_0(u^2)-\kappa^2 /2\Bigr) u_x =0 \;.
$$
This forces $ u_x $ to be identically vanishing on $ ]0,T[ $ which contradicts that $ u(0)= u_0 $ is not a constant and 
 $ u\in C_{w}(\R;L^2(\T)) $. 
 
Note that in the contradiction process, we can replace the assumption on the continuity of the solution-map by the assumption that the flow-map $ u_0\mapsto u(t)$ is continuous from $ L^2(\T) $, equipped with its weak topology, into 
${\mathcal  D}' (\T) $ for all $ t\in ]0,T[ $. Since for each $ t\in \R$, $ u_{n_k}(t)   \rightharpoonup u(t) $ in $ L^2(\T) $ we also get a contradiction. This proves assertion {\it ii'')} of Remark \ref{rere}. 
\subsection{The KdV equation}
First, recall that Miura discovered that the Miura map $ M(u):=u'+u^2 $   maps smooth solutions to the defocusing mKdV equation into smooth solutions to the KdV equation. Actually, it was observed in \cite{CKSTT3}Ê that the Riccati map  
$$
R(u):=u'+u^2-P_0(u^2) 
$$
maps smooth solutions to the defocusing version of equation  \eqref{mKdV*} (i.e. with the  $+$ sign in front of the nonlinear term) into smooth solutions to the KdV equation \eqref{KdV}, i.e. if $ u $ is a smooth solution of the defocusing \eqref{mKdV*}  then $ R(u) $ is a smooth solution of  \eqref{KdV}.  Moreover, according to  \cite{KT0},  this map enjoys the following property :
\begin{theorem}[\cite{KT0}] \label{LemmeMiura}
The Riccati map $ R $ is an isomorphism from   $ L^2_0(\T) $ into $ H^{-1}_0(\T) $. 
\end{theorem}
Actually, we will only need  the Riccati map to be a bijection from $ H^\infty_0(\T) $ into itself. For sake of completeness, we give in the appendix the outline  of the proof of this property.\vspace{2mm} 

Now, let $ w_0\in H^\infty_0(\T) $ and $ \kappa\in \R^* $ be given. We  set $ \theta_0:=w_0-\kappa^2 \cos x $,  $ u_0:=R^{-1}(\theta_0)\in H^\infty_0(\T)$\footnote{It is easy to check that $ u_0\in L^2(\T) $ and $ R(u_0)\in H^\infty_0(\T) $ ensure that $ u_0\in H^\infty_0(\T)$.} and 
$$
u_{0,n}:=u_0+\kappa [\cos(nx)+\cos((n+1)x)]\; .
$$
Clearly, $ u_{0,n} \rightharpoonup u_0 $ in $ L^2(\T) $,  $ \|u_{0,n}\|^2_{L^2} \to  \|u_{0}\|^2_{L^2} +2\kappa^2 \pi$ and
 $ R(u_{0,n}) \rightharpoonup \theta_0+\kappa^2 \cos(x)=w_0 $ in $H^{-1}(\T) $. According to Proposition \ref{propo} there exists a subsequence $ \{u_{n_k}\} $ of  the emanating solutions $ \{u_n\} $
 to \eqref{mKdV*}   and $ u\in C_{w}(\R;L^2_0(\T)) $  such that $ u_{n_k}(t)   \rightharpoonup u(t) $ in $ L^2(\T) $ for all $ t\in \R $ and $ u $ is a weak solution to  \eqref{mKdV*}. To identify $ R(u)$ we will need the following lemma 
 \begin{lem} The operator $   u \mapsto u^2-P_0(u^2)  $   is continuous from  $ F^{-1/16,7/16}_{T}$ 
 into $ {\mathcal D}'(]-T,T[\times  \T) $.
\end{lem} 
\proof 
We set 
 $$
 C(u,u):=u^2 -P_0(u^2)=\sum_{\xi\in \Z^*} \Bigl[ \sum_{(\xi_1,\xi_2)\in \Z^2 \atop 
 \xi_1+\xi_2=\xi} \hat{u}(\xi_1) \hat{u}(\xi_2) \Bigr] e^{i\xi x}
 $$
 Taking $ w\in {\mathcal D} (]-T,T[\times\T) $ and extensions $ {\tilde u_i}\in
F^{-1/16,7/16} $ of $ u_i \in F^{-1/16,7/16}_{T} $   such that $
\|{\tilde u_i}\|_{F^{-1/16,7/16}}
 \le 2 \|u_i\|_{F^{-1/16,7/16}_T} $,
 it holds
\begin{eqnarray*}
J& := & \Bigl| \Bigr(w, C(u_1,u_2)\Bigl)_{L^2(\R\times \T)} \Bigr| \\
& \lesssim & \sum_{(k,k_1,k_2)\in \N^3 } \Bigl|  \Bigr(\Delta_k w, C({ \Delta}_{k_1} {\tilde u_1},{ \Delta}_{k_2} {\tilde u_2})\Bigl)_{L^2(\R\times\T)} \Bigr| \\
\end{eqnarray*}
By symmetry we may assume that $ k_1\ge k_2 $.  We set 
$$ v^{k_1,m}_i: = \gamma(2^{k_1}t -m) {\tilde u_i} , \;  i=1,2,  \; m\in \Z,  $$
where  $ \gamma \, :\, \R  \mapsto [0,1] $ is a $ C^\infty $-function with compact support in $ [-1,1] $ satisfying 
   $ 
     \gamma\equiv 1 \mbox{ on } [-1/2,1/2] \mbox{ and } \sum_{m\in \Z} \gamma(t-m)^2\equiv 1 \mbox{ on } \R$.
   We then obtain 
 $$
 J\lesssim 
 \sum_{(k,k_2,k_3)\in \N^3\atop k_1\ge k_2} \sum_{|m|\lesssim T 2^{k_1}} \Bigl|  \Bigr(\Delta_k w, C({ \Delta}_{k_1}
 v_1^{k_1,m},{ \Delta}_{k_2} v_2^{k_1,m})\Bigl)_{L^2(\R\times \T)} \Bigr| \; .
 $$
We separate  $ \N^3 $ into two regions.\\
{\bf 1.} The region : $ k_1\le 4k+4 $. Then by \eqref{strichartz} and \eqref{gh1}, we can write
\begin{eqnarray}
J & \lesssim & 
 \sum_{(k,k_2,k_3)\in \N^3 \atop k_1\ge k_2} \sum_{|m|\lesssim T 2^{k_1}}
 \| \Delta_k w \|_{L^2}  \| \Delta_{k_1} v_1^{k_1,m} \|_{L^4} 
  \| \Delta_{k_2} v_2^{k_1,m} \|_{L^4}\nonumber  \\
   & \lesssim & T \sum_{(k,k_2,k_3)\in \N^3 \atop k_1\ge k_2} 2^{k_1}Ê 2^{2k_1/3}
 \| \Delta_k w \|_{L^2}  \prod_{i=1}^2 \sup_{m\in \Z} 2^{-k_1/3}\| \Delta_{k_i} v_i^{k_1,m} \|_{X^{0,1/3}} 
\nonumber \\
& \lesssim & T (\sup_{k\in\N} 2^{5k/3}Ê\| \Delta_k w\|_{L^2})  \prod_{i=1}^2 \sup_{k_i\in\N} \sup_{m\in \Z} \| \Delta_{k_i} v_i^{k_1,m} \|_{X^{-1/3,1/3}} \nonumber
 \\
 & \lesssim & T \|w\|_{L^2(\R;H^2(\T))} \|u_1\|_{F^{-1/3,1/3}_{T}} \|u_2\|_{F^{-1/3,1/3}_{T}} 
 \label{rrr1}
\end{eqnarray}
which is acceptable. \\
{\bf 2.} The region : $ k_1> 4k+4 $. Note that in this region, $ k_2\ge  k_1-2\ge 4k+2$. In this region we  will need the well-known resonance relation 
  \begin{equation}\label{reso}
  |\sigma-\sigma_1-\sigma_2|=|3\xi \xi_1 \xi_2|\gtrsim  2^{2k_1} \; ,
  \end{equation}
  where $\xi=\xi_1+\xi_2 $ and $ \sigma, \sigma_1,\sigma_2 $ are defined by \eqref{defsigma}.
   We subdivide this region into subregions with respect to the maximum of $ (|\sigma|, |\sigma_i|)$. \\
     $ \bullet $ $ |\sigma|=\max(|\sigma|,|\sigma_1|,\sigma_2|)   $. In this subregion, making use of \eqref{reso},  \eqref{strichartz} and \eqref{gh1},  we get
     \begin{eqnarray}
     J & \lesssim & 
 \sum_{(k,k_2,k_3)\in \N^3 \atop k_1\ge k_2, \, k_1>4k+4} \sum_{|m|\lesssim T 2^{k_1}}
 \| \Delta_k w \|_{L^2}  \| \Delta_{k_1} v_1^{k_1,m} \|_{L^4} 
  \| \Delta_{k_2} v_2^{k_1,m} \|_{L^4}\nonumber \\
   & \lesssim & T \sum_{(k,k_2,k_3)\in \N^3 \atop k_1\ge k_2, \, k_1>4k+4} 2^{k_1}Ê 2^{2k_1/3}
 2^{-2k_1} \| \Delta_k w \|_{X^{0,1}}  \prod_{i=1}^2 \sup_{m\in \Z} 2^{-k_1/3}\| \Delta_{k_i} v_i^{k_1,m} \|_{X^{0,1/3}} 
\nonumber \\
 & \lesssim & T  \,\|w\|_{X^{0,1}} \|u_1\|_{F^{-1/3,1/3}_{T}} \|u_2\|_{F^{-1/3,1/3}_{T}} \; .
  \label{rrr2}
\end{eqnarray}
   $ \bullet $ $ |\sigma_1|=\max(|\sigma|,|\sigma_1|,\sigma_2|)   $. In this subregion, defining $ {\tilde \nu}_l ,\, l\ge k_1 $ as in \eqref{defnu} and making use of \eqref{reso}, \eqref{propro1} and \eqref{gh1}-\eqref{gh2},  we get
  \begin{eqnarray}\label{est33}
J & \lesssim  &  \sum_{(k,k_1,k_2)\in \N^3 \atop k_1\ge k_2, \, k_1>4k+4} \sum_{|m|\lesssim T 2^{k_1}} 
\sum_{l\ge 0, l_2\ge k_1 }  \|Ê \nu_l \widehat{\Delta_k w} \ast  \tilde{\nu}_{l_2}
\check{\widehat{v_2^{k_1,m}}}\|_{L^2(|\xi|\gtrsim  2^{k_1})}Ê\| \Delta_{k_1} v_1^{k_1,m}\|_{L^2}\nonumber\\
   & \lesssim   & T    \sum_{(k,k_1,k_2)\in \N^3 \atop k_1\ge k_2,\, k_1>4k+4} \sum_{l\ge 0, l_2\ge k_1 } 
   2^{l/2} 2^{k_1}  (2^{l_2/4}2^{-k_1/4}  +1) \|\Delta_k w \|_{L^2}\nonumberÊ\\
    &    &\hspace*{40mm}   \sup_{m\in\Z} \| \tilde{\nu}_{l_2}  \widehat{\Delta_{k_2} v_2^{k_1,m}}\|_{L^2}\, 
   2^{-7k_1/8}  \sup_{m\in\Z } \| \Delta_{k_1}  v_1^{k_1,m}\|_{X^{0,7/16}} \nonumber\\
    & \lesssim   & T  \sum_{(k,k_1,k_2)\in \N^3 \atop k_1\ge k_2, \, k_1>4k+4} \sum_{l\ge 0 } 
   2^{l/2} 2^{k_1} 2^{-21 k_1/16}   \|\Delta_k w \|_{L^2}Ê\|Ê
  \sup_{m\in\Z }Ê\|   \widehat{\Delta_{k_2} v_2^{k_1,m}}\|_{X^{0,7/16,1}}\nonumber\\
    &    & \hspace*{60mm}
 \sup_{m\in\Z }\| \Delta_{k_1}  v_1^{k_1,m}\|_{X^{0,7/16}} \nonumber \\
  & \lesssim   & T   \|w\|_{X^{0,1}} \prod_{i=1}^2 \| u_i\|_{F^{-1/16,7/16}_{T}} \; .
  \label{rrr3}
         \end{eqnarray}
 $ \bullet $  $ |\sigma_2|=\max(|\sigma|,|\sigma_1|,\sigma_2|)   $. Then we can proceed exactly as in the preceding case by exchanging the role of $ v_1^{k_1,m} $ and $ v_2^{k_1,m} $. Gathering \eqref{rrr1}, \eqref{rrr2} and \eqref{rrr3}, we obtain the desired continuity result. \qed \vspace{2mm} \\
By Corollary  \ref{coro1},  Lemma \ref{equicontinuous} and  possibly dilation arguments as in the proof of Proposition \ref{propo}, we know that the sequence $ \{u_n\} $ is relatively compact in $F^{-1/16,7/16}_{T}$ for any $ T>0$. 
From the above lemma we thus infer that 
$$
R(u_n)  \rightharpoonup R(u)  \mbox{ weak star in  }ÊL^\infty(\R; H^{-1}(\T))\; .
$$
We proceed now by contradiction.  Let us assume that the solution-map is continuous at $ w_0=\theta_0+\kappa^2\cos x  $  from $ H^{-1}(\T) $ equipped with its weak topology into ${\mathcal  D}'(]0,T[\times \T) $. Since $ w_0\in H^\infty_0(\T) $, we deduce from Lemma \ref{LemmeMiura} that $ w_0=R(\varTheta_0) $ for some $ \varTheta_0\in H_0^{\infty}(\T) $ with $ \varTheta_0\neq u_0$. By the continuity of the solution-map at $ w_0 $, $ R(u) $ must be equal on  $]0,T[Ê$ to the solution of KdV emanating from 
 $ w_0$ which is nothing else but $ R(\varTheta ) $ where $\varTheta $ is the smooth solution to \eqref{mKdV*}  emanating from $ \varTheta_0$. This ensures that $ u(t)\in H^\infty_0(\T) $ for all $ t\in ]0,T[ $ and the  injectivity of $ R $ on $ H^\infty_0(\T) $ then forces $ u=\varTheta $ on $ ]0,T[ $. This contradicts that $ u(0)=u_0\neq \varTheta_0=\varTheta(0) $ and both functions are weakly continous with values in $ L^2(\T)$.\vspace{2mm} \\

{\bf Acknowledgements:}   L.M.  was partially supported by the ANR project
 "Equa-Disp".  
 \section{Appendix}
 \subsection{A simplified proof of  \eqref{strichartz}}
 We give below a very simple proof of \eqref{strichartz}.  Recall that this estimate was first established in \cite{Bourgain1993}. To simplify the notations we take $ \lambda=1 $ (this corresponds to  $2\pi $-periodic in space functions)  but it easy to check that the proof is exactly the same for any $ \lambda\ge 1$ (see for instance Section \ref{section7.2} for the slight modifications in the case $ \lambda\ge 1$). 
 By the triangle inequality, we write 
 \begin{eqnarray*}
 \| v \|_{L^4_{tx}}^2 = \| v^2 \|_{L^2} = \| \widehat{v} \star \widehat{v} \|_{L^2} 
& \lesssim  & \sum_{l_1\ge 0,l_2\ge 0}\Bigl\|  ( \beta_{l_1} |\hat{v}|)\star (\beta_{l_2} |\hat{v})| \Bigr\|_{L^2}\\
& \lesssim &
  \sum_{l_1\ge l_2\ge 0} \Bigl\|  ( \beta_{l_1} |\hat{v}|)\star (\beta_{l_2} |\hat{v}|) \Bigr\|_{L^2} \quad .
\end{eqnarray*}
where the $ \beta_l $, $l\in \N$,  are defined in \eqref{defnunu}.
The proof of \eqref{strichartz} will then  follow from the following lemma.
 \begin{lem} \label{lemAA} Let $ u_1 $ and $ u_2 $ be $ L^2(\Z\times\R) $-real valued functions then for any $ (l_1,l_2)\in \N^2 $,
 \begin{equation}\label{az}
\Bigl\|  ( \nu_{l_1} u_1 )\star( \nu_{l_2} u_2)  \Bigr\|_{L^2} \lesssim
 \Bigl(2^{l_1} \wedge 2^{l_2}\Bigr)^{1/2}  \Bigl( 2^{l_1}\vee 2^{l_2}\Bigr)^{1/6}
 \|  \nu_{l_1} u_1 \|_{L^2}\, \|  \nu_{l_2} u_2\|_{L^2} \quad .
 \end{equation}
 \end{lem}
Indeed, with this lemma in hand, rewritting $ l_1 $ as $ l_1=l_2+l $ with $ l\in \N $, we get the following chain of inequalities
\begin{eqnarray*}
\sum_{l_1\ge l_2\ge 0} \Bigl\|   ( \nu_{l_1} |\hat{v}|)\star (\nu_{l_2} |\hat{v}|)\Bigr\|_{L^2} & \lesssim &
\sum_{l\ge 0} \sum_{l_2\ge 0} 2^{l_2/2} 2^{(l_2+l)/6}
 \|\nu_{l_2+l} \widehat{v}  \|_{L^2} \|\nu_{l_2} \widehat{v} \|_{L^2} \\
 & \lesssim & \sum_{l\ge 0} \sum_{l_2\ge 0 } 2^{l_2/3} \|\nu_{l_2} \hat{v} \|_{L^2}   2^{-l/6} 2^{(l_2+l)/3}
  \|\nu_{l_2+l} \widehat{v} \|_{L^2} \\
   & \lesssim & \sum_{l\ge 0} 2^{-l/6} \Bigl(\sum_{l_2} 2^{l_2/3}  \|\nu_{l_2} \hat{v} \|_{L^2}^2\Bigr)^{1/2}
\Bigl(\sum_{l_2\ge 0 }   2^{(l_2+l)/3}
  \|\nu_{l_2+l} \widehat{v} \|_{L^2}^2\Bigr)^{1/2} \\
 & \lesssim  & \|v\|_{X^{0,1/3}}^2 \quad .
\end{eqnarray*}
It thus remains to prove  Lemma \ref{lemAA}. Following  the arguments given  in \cite{Bourgain1993KP} (see also [\cite{ST}, page 460]Êfor more details) 
we can assume that $ u $ and $ w$ are supported in $ \{(\tau,\xi)\in \R\times \Z_+\} $. Let us recall that these arguments are based on the fact that the operator  $ {\bf j} : L^2(\R\times \Z) \to  L^2(\R\times \Z) $, defined by  ${\bf j}(u)(\tau,\xi)=u(-\tau,-\xi) $, is an isometry of $ L^2(\R\times \Z) $ satisfying, for any  real-valued  $ L^1\cap L^2 $-functions $ u_1 $ and $ u_2$,
 $$
 \|u_1\star u_2 \|_{L^2(\R\times \Z)} = \| u_1\star {\bf j} (u_2) \|_{L^2(\R\times \Z)}\; .
 $$
By Cauchy-Schwarz in $(\tau_1,\xi_1) $ we infer that
\begin{eqnarray*}
\| ( \beta_{l_1} u_1) \star (\beta_{l_2} u_2) \|_{L^2}^2 & =& \int_{\tau} \sum_{\xi\in \Z}  \Bigl|\int_{\tau_1} \sum_{\xi_1\in  \Z} (\beta_{l_1} u_1)(\tau_1,\xi_1)\, (\beta_{l_2} u_2) (\tau-\tau_1,\xi-\xi_1)\Bigr|^2 \\
& \lesssim  & \int_{\tau}Ê\sum_{\xi\in \Z } \alpha( \tau,\xi) \int_{\tau_1} \sum_{\xi_1\in  \Z}\Bigl| (\beta_{l_1} u_1) (\tau_1,\xi_1)
\, (\beta_{l_2} u_2) (\tau-\tau_1,\xi-\xi_1)\Bigr|^2 \\
& \lesssim & \sup_{\tau\in \R,\xi\in \Z_+}  \alpha(\tau,\xi) \, \|\beta_{l_1} u_1\|_{L^2}^2 \|\beta_{l_2}  u_2\|_{L^2}^2 \quad ,
\end{eqnarray*}
where
\begin{eqnarray*}
\alpha(\tau,\xi) 
 & \lesssim &  mes \Bigl\{ (\tau_1,\xi_1)\in \R\times  \Z_+  \, \slash \,  \; \xi-\xi_1\in  \Z_+, \\
  & &\hspace*{1cm} \; \langle  \tau_1-\xi_1^3\rangle \sim 2^{l_1} \mbox{ and }
\langle \tau-\tau_1-(\xi-\xi_1)^3\rangle \sim 2^{l_2} \, \Bigr\} \\
 & \lesssim &  (2^{l_1}\wedge 2^{l_2}) \; \#    A(\tau,\xi)\; ,  \end{eqnarray*}
with $$
   A(\tau,\xi):=\{  \xi_1\ge 0 / \; \xi-\xi_1\ge 0  \mbox{ and } \langle \tau-\xi_1^3-(\xi-\xi_1)^3\rangle \lesssim 2^{l_1}\vee 2^{l_2}  \}  \; .
   $$ 
  We    separate two regions. In the region  $ \xi^3 \ge 2^{l_1}\vee 2^{l_2} $, we notice that    $ \partial_y^2 [\tau-y^3-(\xi-y)^3]=-6\xi $ which leads to 
$$
 \# A(\tau,\xi)\lesssim \Bigl( \frac{2^{l_1}\vee 2^{l_2}}{\xi}\Bigr)^{1/2} + 1
 \lesssim \Bigl( 2^{l_1}\vee 2^{l_2}\Bigr)^{1/3} \; .
 $$
In the region $ 0\le \xi^3\le  2^{l_1}\vee 2^{l_2}$, we use that $ 0\le \xi_1\le \xi $ to obtain that 
$$
 \#A(\tau,\xi) \le \# \{ \xi_1, \;  0\le \xi_1^3\le 2^{l_1}\vee 2^{l_2}\} \lesssim \Bigl( 2^{l_1}\vee 2^{l_2}\Bigr)^{1/3}  \; .
$$
This completes the proof of \eqref{az}.\qed
\subsection{Proof of  \eqref{propro1}-\eqref{propro2}}\label{section7.2}
We take $ \lambda\ge 1 $. Let $ A\subset  \lambda^{-1}\Z $ and let, for any $ \xi\in   \lambda^{-1}\Z $, $ B(\xi)\subset \lambda^{-1} \Z $. To prove  \eqref{propro1}-\eqref{propro2} we first notice that, by Cauchy-Schwarz,
\begin{align*}
&\Bigl\|\int_{\R}\frac{1}{\lambda} \sum_{\xi_1\in B(\xi)} u_1(\tau_1,\xi_1) \, u_2(\tau-\tau_1,\xi-\xi_1) \, d\tau_1 \Bigr\|_{L^2(A\times \R)}^2 \\
& \lesssim   \int_{\tau}Ê\frac{1}{\lambda}\sum_{\xi\in A } \alpha( \tau,\xi) \int_{\tau_1}Ê\frac{1}{\lambda} \sum_{\xi_1\in  B(\xi)}\Bigl| u_1 (\tau_1,\xi_1)
\, u_2  (\tau-\tau_1,\xi-\xi_1)\Bigr|^2 \, d\tau_1\, d\tau\\
& \lesssim  \sup_{\tau\in \R,\xi\in A} \alpha(\tau,\xi) \, \|u_1\|_{L^2}^2 \|u_2\|_{L^2}^2 \quad ,
\end{align*}
where
$$
\alpha(\tau,\xi) 
  \lesssim  mes \Bigl\{ (\tau_1,\xi_1)\in \R\times B(\xi)  \, \slash \,  (\tau_1,\xi_1)\in \supp(u_1) \mbox{ and }
(\tau-\tau_1,\xi-\xi_1) \in \supp(u_2)\Bigr\} \; .
$$
 Therefore, assuming that  $ \supp(u_i) \subset \{(\tau,\xi)/ \langle \tau-\xi^3\rangle \lesssim L_i \}$, we get 
 $$
 \alpha(\tau,\xi) 
  \lesssim (L_1\wedge L_2)  \, mes[   C(\tau,\xi)]\; 
  $$
with $$
   C(\tau,\xi):=\{  \xi_1\in B(\xi) /\langle \tau-\xi_1^3-(\xi-\xi_1)^3\rangle \lesssim L_1\vee L_2 \}  \; .
   $$ 
   \eqref{propro1} then follows by noticing that  $ \partial^2_y [ \tau-y^3-(\xi-y)^3]=-6\xi $ and thus, for any $|\xi|\ge N>0 $, 
   $$
   mes [   C(\tau,\xi)]=\frac{1}{\lambda}    \#    C(\tau,\xi)\lesssim \frac{1}{\lambda}\Bigl[   \Bigl( \lambda \, \frac{(L_1\vee L_2)^{1/2} }{N^{1/2}}\Bigr) \vee 1\Bigr]  \; .
  $$
Finally, \eqref{propro2} follows by noticing that $$ \partial_y [ \tau-y^3-(\xi-y)^3]=-3\Bigl(y+(\xi-y)\Bigr)\Bigl(y-(\xi-y)\Bigr)$$ and thus, on 
 $ B(\xi)=\{\xi_1 \in\lambda^{-1}\Z\, /\,  \Bigl| |\xi_1|-|\xi-\xi_1|\Bigr|\ge N>0\} $, it holds  $$\Bigl|\Bigl(\xi_1+(\xi-\xi_1)\Bigr)\Bigl(\xi_1-(\xi-\xi_1)\Bigr)\Bigr| \gtrsim N^2\; $$
 which  leads to 
$$
   mes[C(\tau,\xi)]\lesssim \frac{1}{\lambda} \Bigl[ \Bigl(\lambda\,  \frac{L_1\vee L_2 }{N^2}\Bigr)\vee 1\Bigr] \; .
  $$
 \subsection{Outline of the proof of the bijectivity of  the Riccati map from $ H^\infty_0(\T) $ into itself.}
We follow the arguments in \cite{CKSTT3}. For  $ u\in H^\infty_0(\T) $ we denote by $ L_u $ the Schr\"odinger operator with potential $ u$, i.e. 
 $$
 L_u:= -\frac{d}{dx^2} + u 
 $$
 with domain $ H^2(\T) $. 
 One can associate to $L_u $ the energy  $ E_u(\cdot) $ defined on $ H^1(\T) $ by 
 $$
 E_u(\phi) := \langle L_u \phi,\phi\rangle= \int_{\T} \phi_x^2+u \phi^2 \; .
 $$ 
 Since $ L_u $ is a self adjoint operator with compact resolvent, it has a discrete spectrum $ \lambda_1\le \lambda_2\le ... $ with $ \lambda_n \to +\infty $. By the definition of $ E_u $ one must have $ E_u(\phi) \ge \lambda_1 \int_{\T} \phi^2 $ for any $ \phi\in H^1(\T) $ with equality if and only if $ \phi $ is a $ \lambda_1 $-eigenfunction. For $ u\not \equiv 0$, noticing that $ E(1)=0 $ and that 
  $ 1$ is not an eigenfunction, it follows that the first eigenvalue $ \lambda_1 $ is negative. Then, by standard arguments, one can check that $ \lambda_1 $ is a simple eigenvalue with an  eigenfunction 
  that is a non vanishing   $H^\infty(\T)$-function. On the other hand, if $ u\equiv 0 $ it is well-known that the first eigenvalue of $ L_0 $ is $ 0 $ and that the associated eigenspace is spanned by $ 1$. In both cases,  we normalized this eigenfunction by requiring it to be positive and $ L^2$-normalized and we call it $ \phi_1$. Introducing the logarithmic derivative $ v$ of $ \phi_1 $, defined by 
  $$
  v:= \frac{d}{dx} \ln(\phi_1) = \frac{\phi_{1,x}}{\phi_1}\in H^\infty_0(\T) 
  $$
  we observe that 
  $$
  v_x=  \frac{\phi_{1,xx}}{\phi_1}-\Bigl( \frac{\phi_{1,x}}{\phi_1} \Bigr)^2=u-\lambda_1-v^2
  $$
  and thus $ u=v_x+v^2+\lambda_1 $. Taking the means of both sides of this equality it leads to $ \lambda_1=-P_0(v^2) $ which ensures that 
   $ u=R(v) $. This proves that $ R $ is surjective from  $ H^\infty_0(\T) $ into itself. Now, let $ w\in H^\infty_0(\T) $ be such that $ R(w)=u $. Setting  
    $$
    \rho:= \exp(\int_0^x w(s) \, ds) 
    $$
     it is easy to check that $ w=\rho'/\rho $. Observing that 
     $$
     (\frac{d}{dx}+w)(-\frac{d}{dx}+w)=-\frac{d^2}{dx^2} + R(w)+P_0(w^2)=-\frac{d^2}{dx^2} +u+P_0(w^2)\, ,
     $$
     easy calculations then  lead to 
     $$
     E_{u}(\phi)=\int_{\T} (\phi'-w\phi)^2-P_0(w^2) \int_{\T} \phi^2,\quad \quad  \forall \phi\in H^1(\T)  \; .
    $$
    It follows that $ E(\rho)=-P_0(w^2)  \int_{\T} \rho^2 $ and $E(\phi_1)\ge -P_0(w^2)  \int_{\T} \phi_1^2 $. This ensures that 
     $ -P_0(w^2)=\lambda_1 $. Therefore, 
   $\rho/\sqrt{\int \rho^2} =\phi_1 $ and thus $ w=\phi_1'/\phi_1=v $. This proves the injectivity of $ R $ in $ H^\infty_0(\T) $.

  \end{document}